\documentclass[11pt]{article}
\usepackage{amsthm,amsmath,latexsym,amssymb}

\newcommand{\bP}{{\rm |\kern-.15em P}}
\newcommand{\Q}{\kern.3em\rule{.07em}{.65em}\kern-.3em{\rm Q}}
\newcommand{\R}{{\rm I\kern-.15em R}}
\newcommand{\D}{{\rm |\kern-.15em D}}
\newcommand{\h}{{\rm |\kern-.15em H}}
\newcommand{\C}{\kern.3em\rule{.07em}{.65em}\kern-.3em{\rm C}}
\newcommand{\T}{{\rm T\kern-.35em T}}

\theoremstyle{plain}
\newtheorem{theorem}{Theorem}[section]

\newtheorem{proposition}[theorem]{Proposition}
\newtheorem{corollary}[theorem]{Corollary}

\theoremstyle{definition}
\newtheorem{definition}[theorem]{Definition}
\newtheorem{example}[theorem]{Example}

\newtheorem{problem}[theorem]{Problem}

\theoremstyle{remark}
\newtheorem{remark}[theorem]{Remark}

\begin{document}
\title{Composition arithmetic for locally one-to-one entire mappings}
\author{Ronen Peretz}
 
\maketitle

\begin{abstract}
This paper describes a part of the factorization theory of the family of all the entire functions with non vanishing
derivatives. In particular it proves that this family of mappings contains primes. This assures that this family of
entire functions has two non degenerate fractal representations.
\end{abstract}

\section{Introduction}
We consider the class of all the entire functions $f\,:\,\mathbb{C}\rightarrow\mathbb{C}$ with a non vanishing derivative 
and normalized by $f^{'}(0)=1$. We denote this class by ${\rm elh}(\mathbb{C})$ (the class of entire local homeomorphisms).
It contains all the normalized holomorphic $\mathbb{C}$-automorphisms, $z+a$, $a\in\mathbb{C}$. This subclass of 
${\rm elh}(\mathbb{C})$ we denote by ${\rm Aut}(\mathbb{C})$. The set of all the prime entire functions in ${\rm elh}(\mathbb{C})$
is denoted by $\mathcal{M}$. Here the notion of primeness is with respect to composition of mappings as the binary operation. \\
The main result (from our view point) is that $\mathcal{M}\ne\emptyset$ (Theorem \ref{thm8}). \\
In fact, all the functions of the form $f(z)=\int^{z}e^{p(t^{2})}dt$, where $p(t)\in\mathbb{C}[t]-\mathbb{C}$, belong to
${\rm elh}(\mathbb{C})-{\rm Aut}(\mathbb{C})$, satisfy $f(\mathbb{C})=\mathbb{C}$ (using a beautiful idea of E. Calabi) and are
prime entire mappings (Corollary \ref{cor1}). \\
We also prove that a function of the form $f(z)=\int^{z}e^{p(t)}dt$, where $p(t)\in\mathbb{C}[t]-\mathbb{C}$, is a prime entire
mapping, if and only if $f(\mathbb{C})=\mathbb{C}$ (Theorem \ref{thm7}). \\
Moreover, a function of such a form $f(z)=\int^{z}e^{p(t)}dt$ has no non degenerate factorization of the type $f(z)=g(L(z))$, where
$L(z),\,g(z)\in {\rm elh}(\mathbb{C})-{\rm Aut}(\mathbb{C})$ and where $L$ has no Picard exceptional value, i.e. 
$L(\mathbb{C})=\mathbb{C}$ (Theorem \ref{thm6}). \\
Another general result we prove is that if $f\in {\rm elh}(\mathbb{C})$ has a Picard exceptional value, then $f$ is not a prime 
entire function (Theorem \ref{thm4}). \\
In fact such a function has a factorization that includes the logarithm function of a non vanishing entire function which
by itself belongs to ${\rm elh}(\mathbb{C})$. In the special case in which $f(z)=\int^{z}e^{p(t)}dt$ ($p(t)\in
\mathbb{C}[t]-\mathbb{C}$) has a Picard exceptional value, $f(z)$ has one more factorization. It uses the fractional power
$f(z)^{1/(N+1)}$ for some $N\in\mathbb{Z}^{+}$ (See the proof of Theorem \ref{thm7}). \\
Our tools include, among other things, certain systems of initial value problems of some functional-differential equations.
As differential equations they are non linear and of order $1$. Those systems have entire solutions if and only if
the given mapping $f(z)\in  {\rm elh}(\mathbb{C})$ is factorisable (i.e. non prime). In fact, if $f(z)=\int^{z}e^{F(t)}dt\in
{\rm elh}(\mathbb{C})$, where $F(z)\not\equiv {\rm Const.}$ is entire, then we inquire when is $f(z)=g(L(z))$, where
$L(z)\in {\rm elh}(\mathbb{C})-{\rm Aut}(\mathbb{C})$, and where $g(z)$ is an entire non $\mathbb{C}$-automorphism. The
dynamical systems of initial value problems that control the factorization of $f(z)$ are of two possible types: \\
\underline{\bf Type 1:} If $L(\mathbb{C})=\mathbb{C}$ (no Picard exceptional value), then $g\in {\rm elh}(\mathbb{C})$
and hence $g(z)=\int^{z}e^{G(t)}dt$ for some $G(t)\not\equiv {\rm Const.}$ entire. The system is:
$$
\left\{\begin{array}{l} L^{'}(z)=e^{F(z)}\cdot e^{-G(L(z))} \\ L(\alpha)=a \end{array}\right..
$$
\underline{\bf Type 2:} If $L(\mathbb{C})=\mathbb{C}-\{ a\}$ ($a$ is a Picard exceptional value), then 
$g(z)=\int^{z}e^{G(t)}\cdot (t-a)^{N}dt$ for some $N\in\mathbb{Z}^{+}\cup\{ 0\}$. The system is:
$$
\left\{\begin{array}{l} L^{'}(z)=e^{F(z)}\cdot (L(z)-a)^{-N}\cdot e^{-G(L(z))} \\ L(\alpha)=a \end{array}\right..
$$
All this is explained in Section 4 and in Section 3. \\
Thus an important question is to understand when those two dynamical systems have entire solutions $G$ and $L$. Our
approach is to fix arbitrary entire $G(z)$ and analyze the structure of the singular locus of $L(z)$. This locus
is empty if and only if $L(z)$ is also entire, in which case we have our factorization over ${\rm elh}(\mathbb{C})$. \\
The geometric and the topological structures of ${\rm sing}(L)$, are given in Theorem \ref{thm2} and in Theorem \ref{thm3}.
We make a novel use of the generating initial value systems above by varying the data $\alpha$ and $a$ in an appropriate
manner, in order to construct the Riemann surface that uniformizes $L(z)$. These results and technique are of independent
interest of the factorization problem. \\
In fact our motivation for this research originated in the two fractal representations (right an left) of ${\rm elh}(\mathbb{C})$
(Theorem \ref{thm1}):
$$
{\rm elh}(\mathbb{C})=W_{R}\cup\bigcup_{p\in\mathcal{M}}R_{p}({\rm elh}(\mathbb{C})),
$$
$$
{\rm elh}(\mathbb{C})=W_{L}\cup\bigcup_{p\in\mathcal{M}}L_{p}({\rm elh}(\mathbb{C})).
$$
The main result of this paper from our view point, $\mathcal{M}\ne\emptyset$ (Theorem \ref{thm8}) assures that the fractal
representations of ${\rm elh}(\mathbb{C})$ do not degenerate to ${\rm elh}(\mathbb{C})=W_{R}$ or ${\rm elh}(\mathbb{C})=W_{L}$.

Our methods do not seem to work for functions $f(z)=\int^{z}e^{F(t)}dt$ where $F(t)\not\in\mathbb{C}[t]-\mathbb{C}$. Here is
a typical problem left open: \\
Consider an example of Rolf Nevanlinna $f(z)=\int_{0}^{z}e^{e^{t}}dt$. Nevanlinna shows that $f(\mathbb{C})=\mathbb{C}$. Is
that $f(z)$ a prime? In other words can we extend Theorem \ref{thm7} to non polynomial exponents?

\section{Definitions and back ground}

\begin{definition}\label{def1}
We start by recalling few definitions and notations that were used in \cite{r}. \\
${\rm elh}(\mathbb{C})=\{f:\,\mathbb{C}\rightarrow\mathbb{C}\,:\,f\,{\rm is}\,{\rm entire},\,\forall\,z\in\mathbb{C},\,
f^{'}(z)\ne 0,\,f^{'}(0)=1\}$, \\
\\
${\rm Aut}(\mathbb{C})=\{z+a\,:\,a\in\mathbb{C}\}$. \\
\\
For $f\in {\rm elh}(\mathbb{C})$, $A(f)=\{\alpha\in\mathbb{C}\,|\,\exists\,\gamma:\,[0,\infty)\rightarrow\mathbb{C}\,
{\rm continuous},\,{\rm such}\,{\rm that:}\,\lim_{t\rightarrow\infty}|\gamma(t)|=\infty\,\wedge\,
\lim_{t\rightarrow\infty}f(\gamma(t))=\alpha\}$. The curve $\gamma$ is called an asymptotic tract of the asymptotic value $\alpha$. 
$A(f)$ is called the asymptotic variety of the mapping $f\in {\rm elh}(\mathbb{C})$. \\
\\
$T_{f}(f)=\{\gamma:\,[0,\infty)\rightarrow\mathbb{C}\,|\,\gamma\,{\rm is}\,{\rm continuous},\,\lim{t\rightarrow\infty}|\gamma(t)|=\infty\,
\wedge\,\lim_{t\rightarrow\infty}f(\gamma(t))\in\mathbb{C}\}$. $T_{f}(f)$ is the set of all the asymptotic tracts of finite asymptotic
values of the mapping $f\in {\rm elh}(\mathbb{C})$. \\
\\
Two asymptotic tracts $\gamma,\beta\in T_{f}(f)$ of $f\in {\rm elh}(\mathbb{C})$ are equivalent if: \\
a. $\lim_{t\rightarrow\infty}f(\gamma(t))=\lim_{t\rightarrow\infty}f(\beta(t))$. \\
b. $\gamma$ and $\beta$ are homotopic via an homotopy that fixes $\infty$. Thus $\exists\,
H(s,t):\,[0,1]\times [0,\infty)\rightarrow\mathbb{C}$ continuous in $(s,t)$ and satisfying $H(0,t)=\gamma(t)$,
$H(1,t)=\beta(t)$, $\forall\,t\ge 0$, and $\lim_{t\rightarrow\infty}|H(s,t)|=\infty$, $\lim_{t\rightarrow\infty}f(H(s,t))=
\lim_{t\rightarrow\infty}f(\gamma(t))$, $\forall\,0\le s\le 1$. This is an equivalence relation on $T_{f}(f)$. \\
\\
$T_{0f}(f)$ will denote the set of all the equivalence classes of $T_{f}(f)$ modulo the above equivalence relation. \\
\\
For $f\in {\rm elh}(\mathbb{C})$, the right shift mapping on ${\rm elh}(\mathbb{C})$ induced by $f$ is: \\
$R_{f}:\,{\rm elh}(\mathbb{C})\rightarrow {\rm elh}(\mathbb{C}),\,\,\,R_{f}(g)=g\circ f$. \\
\\
The left shift mapping on ${\rm elh}(\mathbb{C})$, induced by $f$ is: \\
$L_{f}:\,{\rm elh}(\mathbb{C})\rightarrow {\rm elh}(\mathbb{C}),\,\,\,L_{f}(g)=f\circ g$.
\end{definition}
\noindent
{\bf Results from \cite{r}:} \\
1. If $f,g\in {\rm elh}(\mathbb{C})$ then $T_{f}(g)\subseteq T_{f}(f\circ g)$, $f(A(g))\subseteq A(f\circ g)$. \\
\\
2. If $f,g\in {\rm elh}(\mathbb{C})$, then $\gamma\in T_{f}(f\circ g)$ implies that either $\gamma\in T_{f}(g)$ or
else $\lim_{t\rightarrow\infty}|g(\gamma(t))|=\infty$. \\
\\
3. If $f,g\in {\rm elh}(\mathbb{C})$, then $A(f\circ g)=A(f)\cup f(A(g))$. \\
\\
4. Let $f\in {\rm elh}(\mathbb{C})$. If $\exists\,g\in {\rm elh}(\mathbb{C})$ such that $T_{f}(g)=T_{f}(f\circ g)$, then
$f(\mathbb{C})=\mathbb{C}$. \\
\\
5. The following are equivalent: (a) $f\not\in {\rm Aut}(\mathbb{C})$. (b) $R_{f}({\rm elh}(\mathbb{C})\subset
{\rm elh}(\mathbb{C})-{\rm Aut}(\mathbb{C})$. (c) $L_{f}({\rm elh}(\mathbb{C}))\subset {\rm elh}(\mathbb{C})-{\rm Aut}(\mathbb{C})$. \\
\\
6. $\forall\,f\in {\rm elh}(\mathbb{C})$, $R_{f}$ is injective. \\
\\
7. $\forall\,f,g,h\in {\rm elh}(\mathbb{C})$, if $L_{f}(g)=L_{f}(h)$ for some $g\ne h$, then there exists an entire function
$\psi(z)$, such that $g(z)=h(z)+e^{\psi(z)}$. \\
\\
8. $L_{\exp(2\pi i z)}$ is not injective. \\

\begin{definition}\label{def2}
Let $(Y,d)$ be a metric space and let $X$ be a topological space. For a compact set $K\subseteq X$, $\epsilon>0$ and $f\in Y^{X}$,
we define $B_{K}(f,\epsilon)=\{g:\,X\rightarrow Y\,|\,d(f(x),g(x))<\epsilon,\,\forall\,x\in K\}$. The sets $B_{K}(f,\epsilon)$
form a basis for a topology on $Y^{X}$, called the topology of compact convergence (on compact subsets of $X$). It is denoted
by $\tau_{cc}$. 
\end{definition}
\noindent
9. (a) The topological space $({\rm elh}(\mathbb{C}),\tau_{cc})$ is a path connected space. (b) For a fixed $f\in {\rm elh}(\mathbb{C})$,
$R_{f}({\rm elh}(\mathbb{C}))$ is a closed subset of $({\rm elh}(\mathbb{C}),\tau_{cc})$. \\
\\
10. $\forall\,f\in {\rm elh}(\mathbb{C})-{\rm Aut}(\mathbb{C})$, we have the identity $\partial R_{f}({\rm elh}(\mathbb{C}))=
R_{f}({\rm elh}(\mathbb{C}))$ in $({\rm elh}(\mathbb{C}),\tau_{cc})$. \\
\\
11. $\forall\,f,g\in {\rm elh}(\mathbb{C})$ we have $L_{f}^{-1}(L_{f}(g))=\{g(z)+k_{j}e^{\psi(z)}\,|\,j=0,\ldots,N,\,k_{0}=0\}$
where $N\in\mathbb{Z}^{+}\cup\{0,\infty\}$. Moreover, we have, $N\in\{0,\infty\}$. \\
\\
12. The set $\{k_{j}e^{\phi(z)}\,|\,j\in I, k_{0}=0\}$ is a cyclic subgroup of $(\mathbb{C},+)$. ($I$ is countable). \\
\\
13. If $f\in\mathbb{C}[z]$, then $L_{f}$ is injective. \\
\\
14. $\forall\,f\in {\rm elh}(\mathbb{C})-\mathbb{C}[z]$, $g\in {\rm elh}(\mathbb{C})$ if $|L_{f}^{-1}(L_{f}(g))|>1$, then
there exist an entire $\phi(z)$ in $z$ and an entire $h(z,w)$ in $(z,w)$ so that $f(g(z)+we^{\phi(z)})=f(g(z))+
e^{h(z,w)}\sin\pi w$. In that case there exists an entire $k(z,w)$ in $(z,w)$ so that:
$$
\frac{\partial h(z,w)}{\partial w}\cdot\sin\pi w+\pi\cos\pi w=e^{k(z,w)}.
$$
15. $\forall\,f\in {\rm elh}(\mathbb{C})-\mathbb{C}[z]$, $g\in {\rm elh}(\mathbb{C})$ if $|L_{f}^{-1}(L_{f}(g))|>1$ then
there exist three entire functions $\phi(z)$, $h(z)$ and $L(w)$ so that $f(g(z)+we^{\phi(z)})=f(g(z))+e^{L(w)+h(z)}\sin\pi w$. Also
we have: \\
a.  $L^{'}(w)\sin\pi w+\pi\cos\pi w$ never vanishes. \\
b. The function $-h^{'}(z)f(g(z)+we^{\phi(z)})+(g^{'}(z)+w\phi^{'}(z)e^{\phi(z)})f^{'}(g(z)+we^{\phi(z)})$ is
independent of $w$. \\
c. $-h^{'}(z)f(g(z))+g^{'}(z)f^{'}(g(z))\equiv g^{'}(z)-g(z)\phi^{'}(z)$. \\
d. $(\phi^{'}(z)-h^{'}(z))f^{'}(g(z)+we^{\phi(z)})+(g^{'}(z)+w\phi^{'}(z)e^{\phi(z)})f^{''}(g(z)+we^{\phi(z)})\equiv 0$. \\
\\
16. Let $f\in {\rm elh}(\mathbb{C})$. Then $L_{f}$ is not injective if and only if $f(z)=(1/b)e^{bz}+a$ for $a\in\mathbb{C}$
and $b\in\mathbb{C}^{\times}$. \\
\\
This is Theorem 3.33 in \cite{r}. {\bf ERRATUM:} This theorem can not be true and needs a correction of the following form: \\
$L_{f}$ is note injective if and only if $f(z)=s(z)\circ ((1/b)e^{bz}+a)$ where $s\in {\rm elh}(\mathbb{C})$, $a\in\mathbb{C}$
and $b\in\mathbb{C}^{\times}$. \\
The reason is the following simple principle.

\begin{proposition}\label{prop1}
Let $X$ be any set, and let $(G,\circ)$ be a semi-group of mappings $X\rightarrow X$. The binary operation $\circ$ is
composition of mappings. Let $N(G)=\{f:\,X\rightarrow X\,|\,f\,{\rm is}\,{\rm non-injective}\,f\in G\}$. Then if
$f\in N(G)$, then the $f$-right orbit $G\circ f$, satisfies $G\circ f\subseteq N(G)$.
\end{proposition}
\noindent
{\bf Proof.} \\
By the definition of $N(G)$ and by the assumption that $f\in N(G)$ it follows that there are $x_{1},x_{2}\in X$,
$x_{1}\ne x_{2}$ such that $f(x_{1})=f(x_{2})$. Let $g\in G$ and $h=g\circ f$. Then $h\in G$ (because $G$ is
a semi-group) and $h(x_{1})=(g\circ f)(x_{1})=g(f(x_{1}))=g(f(x_{2}))=(g\circ f)(x_{2})=h(x_{2})$. Hence $h=g\circ f\in N(G)$.
Hence $G\circ f\subseteq N(G)$. $\qed $ \\
\\
The parameters that are relevant to Theorem 3.33 are: $X={\rm elh}(\mathbb{C})$, $G=\{L_{g}\,|\,g\in {\rm elh}(\mathbb{C})\}$
and $L_{f}\in N(G)$. By Proposition \ref{prop1} we have $G\circ L_{f}\subseteq N(G)$. We note that $L_{g}\circ L_{f}=L_{g\circ f}$. Hence
$\forall\,g\in {\rm elh}(\mathbb{C})$, $L_{g\circ ((1/b)e^{bz}+a)}$ is not injective. We note that if $IN(G)=\{
f:\,X\rightarrow X\,|\,f\,{\rm is}\,{\rm injective}\,f\in G\}$ then: \\
a. Both $(IN(G),\circ)$ and $(N(G),\circ)$ are sub-semi-groups of $(G,\circ)$. \\
b. $IN(G)\cap N(G)=\emptyset$, $IN(G)\cup N(G)=G$. \\
c. $G\circ N(G)=\{g\circ f\,|\,g\in G,\,f\in N(G)\}\subseteq N(G)$ and if ${\rm id}_{X}\in G$, then $G\circ N(G)=N(G)$. \\
\\
The list of the 16 results quoted above is from the paper \cite{r}. We now would like to apply the ideas that
were introduced in the preprint \cite{r1} and prior to that in the paper \cite{r2}.
In \cite{r2,r1} the semi-group ${\rm et}(\mathbb{C}^{2})$ of \'etale
polynomial mappings $F:\,\mathbb{C}^{2}\rightarrow\mathbb{C}^{2}\in\mathbb{C}[X,Y]^{2}$ such that $\det J_{F}(X,Y)\equiv 1$
takes the place of the semi-group ${\rm elh}(\mathbb{C})$ in the current manuscript. The group ${\rm Aut}(\mathbb{C}^{2})$
of polynomial automorphisms of $\mathbb{C}^{2}$ takes the place of ${\rm Aut}(\mathbb{C})$ in the current manuscript.

\begin{proposition}\label{prop2}
$\forall\,f\in {\rm elh}(\mathbb{C})$, $\forall\,n\in\mathbb{Z}^{+}\cup\{0\}$ we have $A(f^{\circ (n+1)})=
A(f)\cup\bigcup_{k=1}^{n}f^{\circ k}(A(f))$.
\end{proposition}
\noindent
{\bf Proof.} \\
For $n=0$ the left side is $A(f)$ and the right side is also $A(f)$. We will use induction on $n$ and assume that
$A(f^{\circ (n+1)})=A(f)\cup\bigcup_{k=1}^{n}f^{\circ k}(A(f))$. Hence we obtain $f(A(f^{\circ (n+1)}))=
f(A(f)\cup\bigcup_{k=1}^{n}f^{\circ k}(A(f)))=\bigcup_{k=1}^{n+1}f^{\circ k}(A(f))$. Hence:
\begin{equation}
\label{eq1}
A(f)\cup f(A(f^{\circ (n+1)}))=A(f)\cup\bigcup_{k=1}^{n+1}f^{\circ k}(A(f)).
\end{equation}
Using result number 3 above (quoted from \cite{r}), we have $A(f\circ g)=A(f)\cup f(A(g))$ and with
$(f,g)=(f,f^{\circ (n+1)})$ we get $A(f\circ f^{\circ (n+1)})=A(f)\cup f(A(f^{\circ (n+1)}))$, or equivalently:
\begin{equation}
\label{eq2}
A(f^{\circ (n+2)})=A(f)\cup f(A(f^{\circ (n+1)})).
\end{equation}
By equations (\ref{eq1}) and (\ref{eq2}) we get $A(f^{\circ (n+2)})=A(f)\cup\bigcup_{k=1}^{n+1} f^{\circ k}(A(f))$.
The inductive argument is completed. $\qed $ \\
\\
We would like to point out that result number 7 above (from \cite{r}) has, in fact a topological origin. For that
we will state it a bit differently and give a corresponding proof.

\begin{proposition}\label{prop3}
Let $f\in {\rm elh}(\mathbb{C})$ be such that $L_{f}$ is not injective. Then $\forall\,g\ne h$, $g,h\in {\rm elh}(\mathbb{C})$
such that $L_{f}(g)=L_{f}(h)$ and $\forall\,z\in\mathbb{C}$ we have $g(z)\ne h(z)$.
\end{proposition}
\noindent
{\bf Proof.} \\
By the assumptions we have $f\circ g=f\circ h$, and there are points $z\in\mathbb{C}$ for which $g(z)\ne h(z)$.
Suppose, in order to get a contradiction, that there is a point $w\in\mathbb{C}$ for which $g(w)=h(w)$. Then we
have two types of points in $\mathbb{C}$. Those $z\in\mathbb{C}$ for which $g(z)\ne h(z)$ and the complementary set,
where both sets are non-empty. Let us denote by $N$ the first subset of $\mathbb{C}$, i.e.: $N=\{z\in\mathbb{C}\,|\,
g(z)\ne h(z)\}$. The subset $N$ of $\mathbb{C}$ is an open subset in the strong topology because $g$ and $h$ are
local homeomorphisms and so if $g(z)\ne h(z)$ then $\exists\, O$, an open neighborhood of $z$ in the strong topology
such that $g(O)\cap h(O)=\emptyset$. So the complementary subset to $N$ is a closed non-empty subset of $\mathbb{C}$.
Let $w\in N^{c}$, be a boundary point of $N^{c}$, i.e. $w\in\partial N^{c}$. Let $w_{n}\in N$ satisfy $\lim w_{n}=w$.
Then $\forall\,n\in\mathbb{Z}^{+}$, $g(w_{n})\ne h(w_{n})$, $f(g(w_{n}))=f(h(w_{n}))$, and $g(w)=h(w)$. This implies 
that in any strong neighborhood of $g(w)=h(w)$ there are different points, say $g(w_{n})\ne h(w_{n})$, 
$\forall\,n\in\mathbb{Z}^{+}$ large enough, so that $f(g(w_{n}))=f(h(w_{n}))$. Hence $f$ is not injective in any
strong neighborhood of $g(w)=h(w)$ and hence $f\not\in  {\rm elh}(\mathbb{C})$, a contradiction. $\qed $ \\

\begin{definition}\label{def3}
Let $f\in {\rm elh}(\mathbb{C})$. The number of irreducible components of $A(f)$ will be denoted by ${\rm comp}(f)$.
\end{definition}

\begin{proposition}\label{prop4}
The following are true: \\
a. $\forall\,f,g\in {\rm elh}(\mathbb{C})$ we have ${\rm comp}(f\circ g)\le {\rm comp}(f)+{\rm comp}(g)$. \\
b. Let $f\in {\rm elh}(\mathbb{C})$, $k\in\mathbb{Z}^{+}$, and $f_{k}={\rm comp}(f^{\circ k})$. Then the following
inequality holds $f_{k}\le {\rm comp}(f)\cdot k$. \\
c. Let $f\in {\rm elh}(\mathbb{C})$, then we have:
$$
|\sum_{k=1}^{\infty} f_{k}z^{k}| \le {\rm comp}(f)\cdot\left(\frac{|z|}{(1-|z|)^{2}}\right),\,\,\,|z|<1.
$$
Thus the function under the absolute value on the left hand side of this inequality, is analytic in the
unit disk, $|z|<1$, and is majorized by a multiple (by ${\rm comp}(f)$) of the Koebe function.
\end{proposition}
\noindent
{\bf Proof.} \\
a. This follows by fact number 3, quoted above from \cite{r}. We have $A(f\circ g)=A(f)\cup f(A(g))$, and by
definition \ref{def3}, of ${\rm comp}(f)$. \\
b. By part a above, and an inductive argument we have $\forall\,k,n\in\mathbb{Z}^{+}$, $f_{k+n}\le f_{k}+f_{n}$. \\
c. This follows by $f_{1}={\rm comp}(f)$, by part b above, and by the fact that the Koebe function, $k(z)$, is given by:
$$
k(z)=\sum_{k=1}^{\infty}kz^{k}.\,\,\,\,\,\,\,\,\,\,\qed
$$

\section{Splitting ${\rm elh}(\mathbb{C})$ in a natural way}
We define a relation on pairs $f,g\in {\rm elh}(\mathbb{C})$ as follows: $f$ and $g$ are $R$-related if and only if
$f\not\in R_{g}({\rm elh}(\mathbb{C}))$ and $g\not\in R_{f}({\rm elh}(\mathbb{C}))$. We will denote the fact that
$f$ and $g$ are so related by writing $f\sim_{R}g$. Clearly $f\sim_{R}g$ implies that $g\sim_{R}f$. We consider
subsets of ${\rm elh}(\mathbb{C})$ that are composed of pairwise $R$-related functions. The family of all such subsets
will be denoted by $\mathcal{F}_{R}$. Thus: $\mathcal{F}_{R}=\{A\subseteq {\rm elh}(\mathbb{C})\,|\,\forall\,f,g\in A,\,
f\sim_{R}g\}$. \\
Similarly, We define a second relation on pairs $f,g\in {\rm elh}(\mathbb{C})$ as follows: $f$ and $g$ are $L$-related if and only if
$f\not\in L_{g}({\rm elh}(\mathbb{C}))$ and $g\not\in L_{f}({\rm elh}(\mathbb{C}))$. We will denote the fact that
$f$ and $g$ are so related by writing $f\sim_{L}g$. Clearly $f\sim_{L}g$ implies that $g\sim_{L}f$. We consider
subsets of ${\rm elh}(\mathbb{C})$ that are composed of pairwise $L$-related functions. The family of all such subsets
will be denoted by $\mathcal{F}_{L}$. Thus: $\mathcal{F}_{L}=\{A\subseteq {\rm elh}(\mathbb{C})\,|\,\forall\,f,g\in A,\,
f\sim_{L}g\}$.

\begin{definition}\label{def4}
A mapping $c\in {\rm elh}(\mathbb{C})$ will be called composite, if there are two entire functions $f,g\not\in {\rm Aut}(\mathbb{C}^{2})$
such that $c=f\circ g$. If $p\in {\rm elh}(\mathbb{C})$ is not composite, then it will be called a prime function. We will identify
two  primes $p$ and $q$, if $\exists\,L,M\in {\rm Aut}(\mathbb{C}^{2})$, such that $p=L\circ q\circ M$, and we will write $p\doteq q$.
\end{definition}
\noindent
We denote: $\mathcal{P}_{R}=\{A\in\mathcal{F}_{R}\,|\,\forall\,g\in A,\,g\,{\rm is}\,{\rm a}\,{\rm prime}\}$ and
$\mathcal{P}_{L}=\{A\in\mathcal{F}_{L}\,|\,\forall\,g\in A,\,g\,{\rm is}\,{\rm a}\,{\rm prime}\}$.

\begin{remark}\label{rem1}
We note that if $p$ and $q$ are primes and if $p\not\doteq q$, then $f\sim_{R}g$ and hence $\{p,q\}\in\mathcal{P}_{R}$.
Also $f\sim_{L}g$ and hence $\{p,q\}\in \mathcal{P}_{L}$.
It is true that the relations $\sim_{R}$ and $\sim_{L}$ are not an equivalence relations (they are clearly not reflexive), 
but the relation $\doteq$ is an equivalence relation. Hence we identify any prime $p$ with its equivalence set in
$\{{\rm elh}(\mathbb{C})/\doteq\}$.
Clearly the equivalence class of a prime $p$ can not contain a prime $q$ such the $q\not\doteq p$. $\mathcal{P}_{R}$
is partially ordered by set inclusion and Zorn's Lemma implies that it contains a maximal element $\mathcal{M}$. Hence 
$\mathcal{M}$ equals exactly the set of all the (equivalence classes of the) primes. Likewise, the maximal element of 
$\mathcal{P}_{L}$ is the same $\mathcal{M}$.
\end{remark}
\noindent
Our hope should have been to have the very simple situation where we have the following presentation of ${\rm elh}(\mathbb{C})$:
$$
{\rm elh}(\mathbb{C})=\bigcup_{p\in\mathcal{M}}R_{p}({\rm elh}(\mathbb{C})),
$$
and where for any pair of primes $p$ and $q$, $p\not\doteq q\Longrightarrow R_{p}({\rm elh}(\mathbb{C}))\cap
R_{q}({\rm elh}(\mathbb{C}))=\emptyset$. This unfortunately is false for at least two pivotal reasons which we now
elaborate. Similarly we would have liked to have the following simple situation:
$$
{\rm elh}(\mathbb{C})=\bigcup_{p\in\mathcal{M}}L_{p}({\rm elh}(\mathbb{C})),
$$
and where for any pair of primes $p$ and $q$, $p\not\doteq q\Longrightarrow L_{p}({\rm elh}(\mathbb{C}))\cap
L_{q}({\rm elh}(\mathbb{C}))=\emptyset$. However unfortunately also this is false for similar two pivotal reasons.

\begin{remark}\label{rem2}
The theory of factorization of entire function (with respect to composition of mappings, as the binary operation) is not as
simple as that for $\mathbb{Z}$. It is wild. We do not have at this point a theorem on the existence and on the uniqueness
of the factorization in ${\rm elh}(\mathbb{C})$. Well known simple examples such as $ze^{z}\circ e^{z}=e^{z}\circ(z+e^{z})$
indicate the wildness of the factorization theory of ${\rm elh}(\mathbb{C})$. True, in the example above we have $ze^{z},
z+e^{z}\not\in {\rm elh}(\mathbb{C})$, but still such examples make it plausible that the factorization theory of
${\rm elh}(\mathbb{C})$ is probably wild. We note, though that $ze^{z}\circ e^{z}=e^{z}\circ(z+e^{z})\in
R_{z+e^{z}}\cap R_{e^{z}}$, while it is not difficult to see that $L_{z+e^{z}}\cap L_{e^{z}}=\emptyset$.
\end{remark}
\noindent
We now list few more questions that come up if we want to proceed with the factorization theory of ${\rm elh}(\mathbb{C})$.
To start with one could have thought that an appropriate definition for the notions of a composite mapping, and a prime mapping
in ${\rm elh}(\mathbb{C})$ should be the following:

\begin{definition}\label{def5}
A mapping $c\in {\rm elh}(\mathbb{C})$ will be called composite, if there are two entire functions 
$f,g\in {\rm elh}(\mathbb{C})-{\rm Aut}(\mathbb{C}^{2})$
such that $c=f\circ g$. If $p\in {\rm elh}(\mathbb{C})$ is not composite, then it will be called a prime function. We will identify
two  primes $p$ and $q$, if $\exists\,L,M\in {\rm Aut}(\mathbb{C}^{2})$, such that $p=L\circ q\circ M$, and we will write $p\doteq q$.
\end{definition}
\noindent
The difference between definition \ref{def5} and our original definition \ref{def4} is that in the original definition we allow
for a composite mapping in ${\rm elh}(\mathbb{C})$ to have a composition of non-automorphisms entire functions, not necessarily
within ${\rm elh}(\mathbb{C})$ only. If these definitions are not equivalent, we somehow feel better with definition \ref{def4}.
One reason is that otherwise we might end up calling a composite entire function, a prime within the ${\rm elh}(\mathbb{C})$
theory. It will be the analog of calling all the numbers in $2\mathbb{Z}$, except for $0,\pm 1$ and $\pm 2$ primes. Thus we
adopt definition \ref{def4}. Having said that we find our selves with basic questions that need to be answered. For example:
are there any prime mappings in ${\rm elh}(\mathbb{C})$? In other words, is there an entire function $f(z)$ which is prime and which
satisfies $\forall\,z\in\mathbb{C},\,f^{'}(z)\ne 0$? In theory, we might have all the mappings in ${\rm elh}(\mathbb{C})$ being
composite but their decompositions involve entire primes outside our semi-group. Another well known result that demonstrates how
wild can be the factorization theory of the entire functions, is the result of Patrick Tuen Wai NG, \cite{t}. It constructs
an entire function which is the the composition of infinitely many prime mappings of the form $c\cdot e^{z}+z$, for carefully
chosen positive numbers $c$. Namely, there exists a sequence of positive real numbers $\{c_{n}\}_{n=1}^{\infty}$ such that
the sequence of functions $F_{n}(z)=(c_{n}e^{z}+z)\circ\ldots\circ(c_{1}e^{z}+z)$ converges uniformly on compact subsets of
$\mathbb{C}$ to an entire function $F(z)$. Furthermore, for each $n\in\mathbb{Z}^{+}$, $F(z)=H_{n}(z)\circ(c_{n}e^{z}+z)\circ
\ldots\circ(c_{1}e^{z}+z)$ for some entire function $H_{n}$. Hence, there is no uniform bound on the number of prime factors
$c_{n}e^{z}+z$ in different decompositions of $F$ through the family of transcendental entire functions. In trying to prove
that ${\rm elh}(\mathbb{C})=\bigcup_{p\in\mathcal{M}}R_{p}({\rm elh}(\mathbb{C}))$, one needs a lemma of the following form:
If $g\in {\rm elh}(\mathbb{C})$ then there exist an entire $m$ and a prime $p\in {\rm elh}(\mathbb{C})$ such that $g=m\circ p$.
We note that the result of Tuen Wai NG does not determine (excludes) if the above type of lemma is valid. Similarly, in trying
to prove that ${\rm elh}(\mathbb{C})=\bigcup_{p\in\mathcal{M}}L_{p}({\rm elh}(\mathbb{C}))$, one needs a lemma of the following
form: If $g\in {\rm elh}(\mathbb{C})$ then there exist an entire $m,p\in {\rm elh}(\mathbb{C})$, $p$ a prime, such that $g=p\circ m$.
Here it might be that a result of the form of that proved by Tuen Wai NG will be relevant. Namely, if we could construct
an entire function of the form $F(z)=\ldots\circ p_{n}(z)\circ\ldots\circ p_{1}(z)$, where the factors $p_{n}(z)$ are
primes in ${\rm elh}(\mathbb{C})$, and where the limiting function $F(z)\in {\rm elh}(\mathbb{C})$ then maybe such a function
were a candidate for one that have no representation of the form $p\circ m$. \\
The following is basically Tuen Wai NG's result, in which the sides left and right are flipped.

\begin{proposition} {\rm (Tuen Wai NG)} \label{prop5}
There exists a sequence of positive real numbers $\{c_{n}\}_{n=1}^{\infty}$ such that the sequence of functions
$F_{n}(z)=(c_{1}e^{z}+z)\circ\ldots\circ (c_{n}e^{z}+z)$ converges uniformly on compact subsets of $\mathbb{C}$,
to an entire function $F(z)$
\end{proposition}
\noindent
{\bf Proof.} \\
We define the positive real numbers $\{c_{n}\}_{n=1}^{\infty}$ inductively. Take $c_{1}=1$ and suppose that
$c_{1},\dots,c_{k}$ were defined. Then we have the following estimate:
$$
\max_{|z|\le k}|F_{k+1}(z)-F_{k}(z)|=F_{k+1}(k)-F_{k}(k)=F_{k}(c_{k+1}e^{k}+k)-F_{k}(k)\rightarrow_{c_{k+1}\rightarrow 0^{+}} 0,
$$
because $G_{k}(w)=F_{k}(we^{k}+k)-F_{k}(k)$ is an entire function of $w$ which satisfies $G_{k}(0)=0$. We choose $c_{k+1}>0$
that satisfies $\max_{|z|\le k}|F_{k+1}(z)-F_{k}(z)|\le 2^{-k}$. It follows that $\{F_{n}\}_{n=1}^{\infty}$ converges
uniformly on each fixed closed disk $|z|\le R$ and so it converges uniformly on compact subsets of $\mathbb{C}$ to a
function $F(z)$ which must be an entire function. Clearly $F(n)>n$ and so $F(z)$ is not a constant function. The rest
of the proof is as in \cite{t}, for compositions that start on the right and extend indefinitely to the left. $\qed $ \\
\\
\begin{definition}\label{def6}
We define two subsets of ${\rm elh}(\mathbb{C})$: \\
$$
W_{R}=\{f\in {\rm elh}(\mathbb{C})\,|\,{\rm if}\,f=p\circ q,\,{\rm where}\,p\,{\rm and}\,q\,{\rm are}\,{\rm entire},
p,q\not\in {\rm Aut}(\mathbb{C}),\,{\rm then}\,q\,{\rm is}\,{\rm not}\,{\rm a}\,{\rm prime}\},
$$
$$
W_{L}=\{f\in {\rm elh}(\mathbb{C})\,|\,{\rm if}\,f=p\circ q,\,{\rm where}\,p\,{\rm and}\,q\,{\rm are}\,{\rm entire},
p,q\not\in {\rm Aut}(\mathbb{C}),\,{\rm then}\,p\,{\rm is}\,{\rm not}\,{\rm a}\,{\rm prime}\}.
$$
\end{definition}

\begin{remark}\label{rem3}
${\rm Aut}(\mathbb{C})\cup\mathcal{M}\subseteq W_{R}$ and ${\rm Aut}(\mathbb{C})\cup\mathcal{M}\subseteq W_{L}$.
Also $\forall\,p\in\mathcal{M}$ we have $W_{R}\cap (R_{p}({\rm elh}(\mathbb{C})-{\rm Aut}(\mathbb{C})))=
W_{L}\cap (L_{p}({\rm elh}(\mathbb{C})-{\rm Aut}(\mathbb{C})))=\emptyset$.
\end{remark}

\begin{theorem}\label{thm1}
We have the following two representations of ${\rm elh}(\mathbb{C})$: 
$$
{\rm elh}(\mathbb{C})=W_{R}\cup\bigcup_{p\in\mathcal{M}}R_{p}({\rm elh}(\mathbb{C})),
$$
and
$$
{\rm elh}(\mathbb{C})=W_{L}\cup\bigcup_{p\in\mathcal{M}}L_{p}({\rm elh}(\mathbb{C})).
$$
\end{theorem}
\noindent
{\bf Proof.} \\
We will prove the first identity by getting a contradiction to the assumption that there is a
$g\in ({\rm elh}(\mathbb{C})-W_{R})-\bigcup_{p\in\mathcal{M}}R_{p}({\rm elh}(\mathbb{C}))$. Since
$g\in {\rm elh}(\mathbb{C})-W_{R}$ it follows by the definition of $W_{R}$ that there is at least one pair
$m$ and $q$ of entire functions, not linear, such that $q$ is a prime, and such that $g=m\circ q$. Thus we
deduce that $g\in R_{q}({\rm elh}(\mathbb{C}))$ where $q\in\mathcal{M}$. This contradicts the assumption
that $g\not\in\bigcup_{p\in\mathcal{M}}R_{p}({\rm elh}(\mathbb{C}))$. This proves the first identity. The
second identity can be proved in a similar manner. $\qed $ \\
\\
The next simple example shows that the subsets $W_{R}$ and $W_{L}$ in the two representations, respectively
are not a redundant component. They must be there.

\begin{proposition}\label{prop6}
$e^{z}\in W_{R}-({\rm Aut}(\mathbb{C})\cup\mathcal{M})$, and also $e^{z}\in W_{L}-({\rm Aut}(\mathbb{C})\cup\mathcal{M})$.
\end{proposition}
\noindent
{\bf Proof.} \\
We note the simple identity $e^{z}=z^{2}\circ z^{3}\circ\ldots\circ z^{n}\circ e^{z/n!}$. It implies that
$e^{z}\not\in\mathcal{M}$. Next let us assume that $e^{z}=m(z)\circ p(z)$ for some entire $m(z)$ and some
prime mapping $p(z)$. Then $e^{z}=p^{'}(z)\cdot m^{'}(p(z))$. So $p^{'}(z)$ never vanishes and hence there
exists an entire function $G(z)$ so that $p^{'}(z)=e^{G(z)}$. By order of growth considerations we 
deduce that $G(z)$ must be linear and hence $p(z)=e^{az+b}$, $a\ne 0$. But we noticed at the beginning of
the proof that $e^{az+b}\not\in\mathcal{M}$, which contradicts the assumption we started with. A similar
argument proves the result for left compositions. $\qed $ \\

\begin{remark}\label{rem4}
One criticism to the simple example given in Proposition \ref{prop6} is that indeed $e^{z}$ is not a prime
but it is a pseudo-prime. Thus it is of interest to find an example which is not even a pseudo-prime that
belongs to $W_{R}-({\rm Aut}(\mathbb{C})\cup\mathcal{M})$ and one that belongs to $W_{L}-({\rm Aut}(\mathbb{C})\cup\mathcal{M})$.
Also we note that Theorem \ref{thm1} does not prevent the degenerate situations where we might be
having: ${\rm elh}(\mathbb{C})=W_{R}$ or ${\rm elh}(\mathbb{C})=W_{L}$. These might be true in the
unfortunate case in which $\mathcal{M}=\emptyset$, i.e. there are no prime mappings within ${\rm elh}(\mathbb{C})$.
Even though such a result is interesting, from our point of view the interesting parts in the two
representations of Theorem \ref{thm1} are $\bigcup_{p\in\mathcal{M}}R_{p}({\rm elh}(\mathbb{C}))$ and
$\bigcup_{p\in\mathcal{M}}L_{p}({\rm elh}(\mathbb{C}))$.
\end{remark}

\begin{remark}\label{rem5}
It is worth comparing the representations we have for the semi-group ${\rm elh}(\mathbb{C})$ (in Theorem \ref{thm1})
to the representations we have in the algebraic setting of the semi-group ${\rm et}(\mathbb{C}^{2})$ in \cite{r2} and in \cite{r1}.
Any $F\in {\rm et}(\mathbb{C}^{2})-{\rm Aut}(\mathbb{C}^{2})$ can always be represented as a composition of at most
$d_{F}$ primes in ${\rm et}(\mathbb{C}^{2})$, where $d_{F}$ is the geometric degree of the mapping $F$. Thus in the
algebraic representations of ${\rm et}(\mathbb{C}^{2})$ the place of the complicated subsets $W_{R}$ and $W_{L}$
is taken by ${\rm Aut}(\mathbb{C}^{2})$. The representations in the algebraic setting of ${\rm et}(\mathbb{C}^{2})$
even tough seem to be driven by factorization into primes are divided naturally into two parts that have geometric
flavor. One part contains the automorphisms (which over $\mathbb{C}^{2}$ coincide with the injective mappings) and the second 
part contains the \'etale mappings which are not injective. The Jacobian Conjecture speculates that this second
part is void. In a sense the parallel of that would be (if true) a theorem that will show that there are no prime
mappings in the semi-group ${\rm elh}(\mathbb{C})$. This question was indeed our motivation in the research 
described in this paper.
\end{remark}

\section{What is the structure of a composite mapping in ${\rm elh}(\mathbb{C})$?}

\begin{proposition}\label{prop7}
The mapping $f(z)\in {\rm elh}(\mathbb{C})$ is a composite mapping if and only if there exist two
entire functions $H(z)$ and $G(z)$ that satisfy: \\
if 
$$
\Im\left\{\int_{0}^{z}e^{H(t)}dt\right\}=\mathbb{C},
$$
then
$$
f(z)=\int_{0}^{\int_{0}^{z}\exp\{H(w)\}dw}\exp\{G(t)\}dt,
$$
and if 
$$
\Im\left\{\int_{0}^{z}e^{H(t)}dt\right\}=\mathbb{C}-\{a\},\,\,{\rm for}\,\,{\rm some}\,\,a\in\mathbb{C},
$$
then
$$
f(z)=\int_{0}^{\int_{0}^{z}\exp\{H(w)\}dw}(t-a)^{N}\exp\{G(t)\}dt.
$$
Here $N\in\mathbb{Z}^{+}\cup\{0\}$, and the entire function $H(z)$ is not a constant function, and if $N=0$ then also
$G(z)$ is not a constant function.
\end{proposition}
\noindent
{\bf Proof.} \\
The function $f(z)$ is a composite function if and only if there exist two entire but not linear functions $g$ and $h$
such that $f=g\circ h$. For $f(z)$ to belong to the semi-group ${\rm elh}(\mathbb{C})$ the following condition must be
fulfilled: $f^{'}(z)$ never vanishes and $f^{'}(0)=1$. We can assume without losing the generality that $f(0)=0$.
Since $f^{'}(z)=g^{'}(h(z))\cdot h^{'}(z)$, it follows that both $h^{'}(z)$ and $g^{'}(h(z))$ never vanish, and also
that $g^{'}(h(0))\cdot h^{'}(0)=1$. We deduce that there exists an entire but not a constant function $H(z)$, such that
$h^{'}(z)=\exp\{H(z)\}$ and so we can assume that $h(z)=\int_{0}^{z}\exp\{H(t)\}dt$ (for assuming that $h(0)=0$ is not
a significant restriction). Also, $g^{'}(w)$ can not have a zero on the image of $h$, i.e. on $h(\mathbb{C})$. By the 
Picard Theorem we either have $h(\mathbb{C})=\mathbb{C}$, or $h(\mathbb{C})=\mathbb{C}-\{a\}$ for some $a\in\mathbb{C}$.
In the first case, where $h(\mathbb{C})=\mathbb{C}$, the entire function $g^{'}(w)$ never vanishes on $\mathbb{C}$ and
so there exists an entire but not a constant function $G(z)$ such that $g^{'}(w)=\exp\{G(w)\}$ and this time we must have
$g(w)=\int_{0}^{w}\exp\{G(t)\}dt$ by $h(0)=0$ and by $f(0)=g(h(0))=0$. In the second case there exists an entire function
$G(z)$ such that $g^{'}(w)=(w-a)^{N}\exp\{G(w)\}$, where $N$ is the order of the zero of $g^{'}(w)$ at $w=a$.
If $N=0$ then necessarily the function $G(z)$ can not be a constant function. If $N>0$ then $G(z)$ can be any entire function.
Finally we substitute the representations of the factors $h(z)$ and $g(z)$ into the equation $f=g\circ h$. In the case
where $h(\mathbb{C})=\mathbb{C}$, we saw that we have:
$$
h(z)=\int_{0}^{z}\exp\{H(t)\}dt,\,\,{\rm and}\,\,g(z)=\int_{0}^{z}\exp\{G(t)\}dt.
$$
In this case we obtain:
$$
f(z)=g(h(z))=\int_{0}^{\int_{0}^{z}\exp\{H(w)\}dw}\exp\{G(t)\}dt.
$$
If, on the other hand we have $h(\mathbb{C})=\mathbb{C}-\{a\}$, for some number $a\in\mathbb{C}$, we saw that we have:
$$
h(z)=\int_{0}^{z}\exp\{H(t)\}dt,\,\,{\rm and}\,\,g(z)=\int_{0}^{z}(t-a)^{N}\exp\{G(t)\}dt\,\,{\rm for}\,\,{\rm some}\,\,a\in\mathbb{C}.
$$
In this case we obtain:
$$
f(z)=g(h(z))=\int_{0}^{\int_{0}^{z}\exp\{H(w)\}dw}(t-a)^{N}\exp\{G(t)\}dt.
$$
The factor $h$ is an entire non-linear function if and only if the entire function, $H$ is non-constant. The factor $g$
is an entire non-linear function if and only if the entire function $(t-a)^{N}\exp\{G(t)\}$ is non-constant and this is
equivalent to: either $N>0$ and $G$ is any entire function, or $N=0$ and $G$ is an entire non-constant function. $\qed $ \\

\section{Factorization within ${\rm elh}(\mathbb{C})$, and a certain functional non-linear ode of order one}

\underline{Case 1:} If $\Im\left\{\int_{0}^{z}e^{H(t)}dt\right\}=\mathbb{C}$, then $f(z)=\int_{0}^{\int_{0}^{z}e^{H(w)}dw}
e^{G(t)}dt$, where $H,G\not\equiv {\rm Const.}$ \\
\\
\underline{Case 2:} If $\Im\left\{\int_{0}^{z}e^{H(t)}dt\right\}=\mathbb{C}-\{a\}$, then $f(z)=\int_{0}^{\int_{0}^{z}
e^{H(w)}dw}(t-a)^{N}\cdot e^{G(t)}dt$, where $H\not\equiv {\rm Const.}$, $N\in\mathbb{Z}^{+}\cup\{0\}$ and in case
$N=0$, then also $G\not\equiv {\rm Const.}$. \\
\\
In both cases both $H$ and $G$ are entire functions. If ${\rm elh}(\mathbb{C})$, then $f(z)=\int_{0}^{z}e^{F(t)}dt$ for
some entire function $F$. Hence we get:
$$
\int_{0}^{z}e^{F(t)}dt=\int_{0}^{\int_{0}^{z}e^{H(w)}dw}e^{G(t)}dt,\,\,\,{\rm in}\,{\rm case}\,1,
$$
and
$$
\int_{0}^{z}e^{F(t)}dt=\int_{0}^{\int_{0}^{z}e^{H(w)}dw}(t-a)^{N}\cdot G(t)dt,\,\,\,{\rm in}\,{\rm case}\,2.
$$
On differentiating we get:
$$
e^{F(z)}=e^{H(z)}\cdot e^{G\left(\int_{0}^{z}e^{H(w)}dw\right)},\,\,\,{\rm in}\,{\rm case}\,1,
$$
and
$$
e^{F(z)}=e^{H(z)}\cdot\left(\int_{0}^{z}e^{H(w)}dw-a\right)^{N}\cdot e^{G\left(\int_{0}^{z}e^{H(w)}dw\right)},\,\,\,{\rm in}\,{\rm case}\,2.
$$
Comparing exponents we get:
$$
F(z)=H(z)+G\left(\int_{0}^{z}e^{H(w)}dw\right)+2\pi i k,\,\,\,{\rm in}\,{\rm case}\,1,
$$
\begin{equation}
\label{eq3} {\rm and}
\end{equation}
$$
F(z)=H(z)+G\left(\int_{0}^{z}e^{H(w)}dw\right)+N\log\left(\int_{0}^{z}e^{H(w)}dw-a\right)+2\pi i k,\,\,\,{\rm in}\,{\rm case}\,2.
$$
The problem we face is the following: Given an entire function $F(z)$ can we find two entire functions $H(z)$ and $G(z)$
or at least prove that they exist? More explicitly, in order to prove that ${\rm elh}(\mathbb{C})$ contains no prime mapping,
we need to show that any $f(z)\in {\rm elh}(\mathbb{C})$ can be factored (non-trivially). This is equivalent to the following:
For any entire function $F(z)$ there exist solutions, $H(z)$ and $G(z)$ which are entire functions, for the equation (\ref{eq3}), and
not trivial (as required in case 1 and in case 2). If we differentiate equation (\ref{eq3}) we get:
$$
F^{'}(z)=H^{'}(z)+e^{H(z)}\cdot G^{'}\left(\int_{0}^{z}e^{H(w)}dw\right),\,\,\,{\rm in}\,{\rm case}\,1,
$$
\begin{equation}
\label{eq4} {\rm and}
\end{equation}
$$
F^{'}(z)=H^{'}(z)+e^{H(z)}\cdot G^{'}\left(\int_{0}^{z}e^{H(w)}dw\right)+\left(\frac{Ne^{H(z)}}{\int_{0}^{z}e^{H(w)}dw-a}\right),
\,\,\,{\rm in}\,{\rm case}\,2.
$$
We simplify now the compositions $G(\int_{0}^{z}e^{H(w)}dw)$ and $G^{'}(\int_{0}^{z}e^{H(w)}dw)$ in equations (\ref{eq3}) and
equations (\ref{eq4}) respectively by introducing the function $L(z)=\int_{0}^{z}e^{H(w)}dw$. Then $L(z)\in {\rm elh}(\mathbb{C})$
and $L^{'}(z)=e^{H(z)}$ or equivalently $H(z)=\log L^{'}(z)$. The equation which corresponds to case 1 becomes:
$$
e^{F(z)}=L^{'}(z)e^{G(L(z))},\,\,\,{\rm instead}\,{\rm of}\,e^{F(f)}=e^{H(z)}\cdot e^{G\left(\int_{0}^{z}e^{H(w)}dw\right)}.
$$
So:
\begin{equation}
\label{eq5}
L^{'}(z)=e^{F(z)}\cdot e^{-G(L(z))}.
\end{equation}
The first case in equation (\ref{eq3}), $F(z)=H(z)+G\left(\int_{0}^{z}e^{H(w)}dw\right)+2\pi i k$, becomes:
\begin{equation}
\label{eq6}
F(z)=\log L^{'}(z)+G(L(z))+2\pi i k.
\end{equation}
If we differentiate this we obtain $F^{'}(z)=L^{''}(z)/L^{'}(z)+L^{'}(z)G^{'}(L(z))$, i.e.:
\begin{equation}
\label{eq7}
L^{''}(z)+(L^{'}(z))^{2}G^{'}(L(z))-L^{'}(z)F^{'}(z)=0.
\end{equation}
So in order to prove that $f(z)\in {\rm elh}(\mathbb{C})$ can be factored, i.e., is not a prime, we need to show that
given the entire function $F(z)$ (where $f^{'}(z)=e^{F(z)}$), there exist two entire functions $G(z)$ and $L(z)$ such
that equation (\ref{eq5}) (or (\ref{eq6}) or (\ref{eq7}) or other variants) is satisfied. The equation is a functional differential
equation and we look for entire solutions $G$ and $L$ for this equation. The problem of finding a local holomorphic solution is
well known. One can find a set of results for the first order ordinary differential equations in the complex domain. We
will use the book \cite{hi} as our reference book. We quote a few results from chapter 2 of \cite{hi}. Existence
and Uniqueness Theorems, pages 40-75. Section 2.2: The fixed point method: We are concerned with the equation $w^{'}=F(z,w)$,
where $(z,w)\rightarrow F(z,w)$ is holomorphic in the dicylinder given by: $D:\,|z-z_{0}|\le a$, $|w-w_{0}|\le b$. It is
required to find a function $z\rightarrow w(z;z_{0},w_{0})$, holomorphic in some disk $|z-z_{0}|<r\le a$, such that
$$
\left\{\begin{array}{l} w^{'}(z;z_{0},w_{0})=F(z,w(z;z_{0},w_{0})) \\ w(z_{0};z_{0},w_{0})=w_{0}.\end{array} \right.
$$
This is done with the aid of Banach fixed-point theorem on contractions on complete metric spaces. $F$ satisfies two
conditions, namely, $|F(z,w)|<M$, and $|F(z,u)-F(z,v)|<K\cdot |u-v|$, for suitably chosen $K$ and $M$ and for
$(z,w)$, $(z,w)$ and $(z,v)$ in $D$. \\
\\
{\bf Theorem 2.2.1} (in \cite{hi}) Under the stated assumptions on $F$, in the disk $D_{0}:\,|z-z_{0}|<r$, where
$$
r<\min\left(a,\frac{b}{M},\frac{1}{K}\right),
$$
the equation $w^{'}=F(z,w)$ has a unique holomorphic solution satisfying the initial value condition $w(z_{0};z_{0},w_{0})=w_{0}$. \\
\\
\\
If we consider equation (\ref{eq5}): $L^{'}(z)=e^{F(z)}e^{-G(L(z))}$ where the entire function $F(z)$ is given and if we fix 
non-constant entire function $G(z)$ (arbitrary and completely to our choice), then in the notations of Theorem 2.2.1 of
\cite{hi}, we have $L(z)=w(z)$ and $F(z,w)=e^{F(z)}e^{-G(w)}$. The function $F(z,w)$ is an entire function and hence the radii
$a$, $b$ of our dicylinder are completely arbitrary positive real numbers. The optimal constants $M$ and $K$ are given by:
$$
M=\max_{D}|e^{F(z)}e^{-G(w)}|=\max_{|z-z_{0}|\le a,\,|w-w_{0}|\le b}e^{\Re F(z)}e^{-\Re G(w)},
$$
$$
K=\min\{k\,|\,|e^{F(z)}e^{-G(u)}-e^{F(z)}e^{-G(v)}|<k\cdot |u-v|,\,\forall\,(z,u),\,(z,v)\in D\}=
$$
$$
=\min\{k\,|\,|e^{-G(u)}-e^{-G(v)}|<k\cdot e^{-\Re F(z)}\cdot |u-v|,\,\forall\,(z,u),\,(z,v)\in D\}.
$$
Since $L(z)=\int_{0}^{z}e^{H(w)}dw$ it follows that $L(0)=0$ and we can choose, for example $(z_{0},w_{0})=(0,0)$.
We note that in general $L(z)$ will not be an entire function. By $H(z)=\log L^{'}(z)$ we see that also $H(z)$ is
in general not entire. If it were, then this would have proved the factorization:
$$
f(z)=\int_{0}^{\int_{0}^{z}e^{H(w)}dw}e^{G(t)}dt=\int_{0}^{L(z)}e^{G(t)}dt.
$$
If we define (as we did in the entire case), $g(w)=\int_{0}^{w}e^{G(t)}dt$ and $h(z)=\int_{0}^{z}e^{H(w)}dw=L(z)$, then we 
obtain the factorization: $f(z)=g(h(z))$, $|z|<r$. This completes the proof of the following:
\begin{proposition}\label{prop8}
Let $f(z)\in {\rm elh}(\mathbb{C})$. Let $g(z)\in {\rm elh}(\mathbb{C})$. Then there exist an $r>0$ and an $h(z)\in
H(\{z\in\mathbb{C}\,|\,|z|<r\})$ such that $f(z)=g(h(z))$, $|z|<r$. The radius $r$ can be any positive real number
smaller than or equal to $\min\left(a,b/M,1/K\right)$, where $a,b>0$ are arbitrary. Once $a$ and $b$ were chosen we 
have:
$$
M=\max_{|z|\le a,\,|w|\le b}\left|\frac{f^{'}(z)}{g^{'}(w)}\right|
$$
and
$$
K=\min\left\{k\,|\,\left|\frac{f^{'}(z)}{g^{'}(u)}-\frac{f^{'}(z)}{g^{'}(v)}\right|<k\cdot |u-v|,\,\forall\,
|z|\le a,\,|u|\le b,\,|v|\le b\right\}.
$$
\end{proposition}

\begin{remark}\label{rem6}
This proposition is not surprising for it is the claim of the implicit function theorem, except that we have here
quantitative estimate for the radius $r$ in which circle the factorization takes place. We note that even though 
$a,b >0$ are arbitrary, we do not want them to be too large because:
$$
\lim_{(a,b)\rightarrow (\infty,\infty)}\frac{b}{M}=\lim_{(a,b)\rightarrow (\infty,\infty)}\frac{1}{K}=0.
$$
\end{remark}
\noindent
Section 2.3 (of \cite{hi}): The method of successive approximations, by \'Emile Picard adds the information that
in fact in Proposition \ref{prop8} we can take $r<\min\left(a,b/M\right)$. The condition $rk<1$ is no longer
imposed. \\
Section 2.4 (of \cite{hi}): Majorants and majorant methods, uses the idea that apparently was due to Issac Newton.
The result is the following: Let $F(z,w)=\sum_{j=0}^{\infty}\sum_{k=0}^{\infty} c_{jk}z^{j}w^{k}$, and let
$G(z,w)=\sum_{j=0}^{\infty}\sum_{k=0}^{\infty}C_{jk}z^{j}w^{k}$ be both analytic in $|z|\le a$, $|w|\le b$.
Also we assume that $G$ is a majorant of $F$, designated by $F(z,w)\ll G(z,w)$ and defined by the infinite 
set of inequalities: $|c_{jk}|\le C_{jk}$, $\forall\,j,k\in\{0,1,2,3,\ldots\}$. Suppose that: 
$$
W^{'}(z)=G(z,W(z)),\,\,\,W(0)=0,
$$
has a solution $W(z)=\sum_{j=1}^{\infty}C_{j}z^{j}$ convergent for $|z|<r$. Let: $w(z)=\sum_{j=1}^{\infty}c_{j}z^{j}$,
be a formal solution of $w^{'}(z)=F(z,w(z))$, $w(0)=0$. Then Theorem 2.4.1 (in \cite{hi}) asserts that if $w(z)\ll W(z)$ and 
the series of $w(z)$ is absolutely convergent for $|z|<r$ and is the unique solution of it's system. \\
Section 2.5 (of \cite{hi}): The Cauchy majorant: is a specific construction of Cauchy that uses the general result
of section 2.4. Also in this section we consider the system: $w^{'}(z)=F(z,w(z))$, $w(0)=0$, and it is proved that
this has a series solution $\sum_{n=1}^{\infty}c_{n}z^{n}$, which is absolutely convergent in the disk $|z|<R$, where
Theorem 2.5.1 (of \cite{hi}) asserts that,
$$
R=q\left(1-\exp\left(-\frac{b}{2aM}\right)\right),
$$
where: $M\equiv M(a,b)=\sum_{j=0}^{\infty}\sum_{k=0}^{\infty}|c_{jk}|a^{j}b^{k}$, with $F(z,w)=\sum_{j=0}^{\infty}\sum_{k=0}^{\infty}
c_{jk}z^{j}w^{k}$. This might tempt to be used as follows: Let $f(z)\in{\rm elh}(\mathbb{C})$ and let $f^{'}(z)=e^{F(z)}$
for some entire function $F(z)$. If we could have found an entire function $G(z)$ satisfying the following two conditions: \\
1. $e^{F(z)}e^{-G(w)}=\sum_{j=0}^{\infty}\sum_{k=0}^{\infty}c_{jk}z^{j}w^{k}$ where $c_{jk}\ge 0$, $\forall\,j,k\in\{0,1,2,3,\ldots\}$. \\
2. $\lim_{b\rightarrow\infty}be^{\Re G(b)}=\infty$. \\
Then $f(z)$ has a factorization $f(z)=g(h(z))$, where $g^{'}(w)=e^{G(w)}$. \\
\\
{\bf Proof.} \\
By condition 1 we have $M(a,b)=e^{F(a)}e^{-G(b)}$ for any $a,b>0$. This is because our differential equation is:
$$
L^{'}(z)=F(z,L(z))=e^{F(z)}e^{-G(L(z))},
$$
where, since $L(z)=\int_{0}^{z}e^{H(w)}dw$, it follows that $L(0)=0$, just like in the system of Theorem 2.4.1 of \cite{hi}. Thus 
the function $F(z,w)=e^{F(z)}e^{-G(w)}$ is a product of two entire functions, hence is itself is an entire function and hence we
can freely choose any $a,b>0$ to define our dicylinder: $D\,:\,|z|\le a$, $|w|\le b$. We define the formal solution
by $\sum_{n=1}^{\infty}c_{n}z^{n}$ and we would like to estimate its radius of convergence. As Cauchy suggests we take:
$$
M\equiv M(a,b)=\sum_{j=0}^{\infty}\sum_{k=0}^{\infty}|c_{jk}|a^{j}b^{k}=\sum_{j=0}^{\infty}\sum_{k=0}^{\infty}c_{jk}a^{j}b^{k}=
e^{F(a)}e^{-G(b)}.
$$
The one before the last equality follows by condition 1. We deduce from Theorem 2.5.1 that the formal solution
converges for:
$$
|z|<R=a\left(1-\exp\left(-\frac{b}{2aM}\right)\right)=a\left(1-\exp\left(-\frac{be^{G(b)}}{2ae^{F(a)}}\right)\right).
$$
We will show that $\exists\,\epsilon>0$, a constant (i.e. independent of $a$ and $b$) such that:
$$
1-\exp\left(-\frac{1}{2}\cdot\frac{be^{G(b)}}{ae^{F(a)}}\right)\ge\epsilon>0,
$$
for appropriate choice of $(a,b)$ for which $(a,b)\rightarrow\infty$. In fact $b=b(a)$ where $\lim_{a\rightarrow\infty}b(a)=\infty$.
It will show that, then $R=\infty$ and hence the formal solution is in fact an entire function $L(z)=\int_{0}^{z}e^{H(w)}dw=
\sum_{n=1}^{\infty}c_{n}z^{n}$. This would prove what we wanted. Thus we need:
$$
\exp\left(-\frac{1}{2}\cdot\frac{be^{G(b)}}{ae^{F(a)}}\right)\le 1-\epsilon,\,\,\,{\rm so}\,
-\frac{1}{2}\cdot\frac{be^{G(b)}}{ae^{F(a)}}\le\log(1-\epsilon).
$$
Thus we need:
$$
\frac{b\cdot e^{G(b)}}{a\cdot e^{F(a)}}\ge -2\log(1-\epsilon).
$$
Let us denote $\delta=-2\log(1-\epsilon)$, then $\delta>0$ is a constant (independent of $a$ and $b$) provided that $0<\epsilon<1$
exists. By condition 2 we have: $\lim_{b\rightarrow\infty}|be^{G(b)}|=\lim_{b\rightarrow\infty} b|G(b)|=\lim_{b\rightarrow\infty}
be^{\Re G(b)}=\infty$. Hence, we can indeed define a function $b=b(a)$, which is positive and increasing and satisfies
for $a\rightarrow\infty$ the equation $b(a)\rightarrow\infty$ and also, say, 
\begin{equation}
\label{eq8}
\frac{be^{G(b)}}{ae^{F(a)}}\ge 1.
\end{equation}
Thus we can take $\delta=1$, i.e. $0<\epsilon=1-e^{-1/2}<1$. The reason for the existence of such a $b=b(a)$ for which
inequality (\ref{eq8}) is satisfied follows by the following reasoning: If $\phi(x)$ and $\psi(x)$ are strictly
increasing functions such that $\lim_{x\rightarrow\infty}\phi(x)=\lim_{x\rightarrow\infty}\psi(x)=\infty$, then there
exists a function $\theta(x)$ which is increasing and satisfies $\lim_{x\rightarrow\infty}\theta(x)=\infty$ and
$\phi(\theta(x))/\psi(x)\ge 1$. We can simply take $\theta(x)=\phi^{-1}(\psi(x))$. $\phi^{-1}$ exists because $\phi$ is
an increasing function. Our functions $\phi(b)=be^{G(b)}$ and $\psi(a)=ae^{F(a)}$ satisfy the requirements. $\qed $

\begin{remark}\label{rem7}
Condition 1 above can be written more explicitly. If $e^{F(z)}=\sum_{j=0}^{\infty}\alpha_{j}z^{j}$ and 
$e^{-G(w)}=\sum_{k=0}^{\infty}\beta_{k}w^{k}$, then the variables $z$ and $w$ are separated and the product of the functions
$e^{F(z)}e^{-G(w)}$ behaves like the standard scalar product in the following sense:
$$
e^{F(z)}e^{-G(w)}=\left(\sum_{j=0}^{\infty}\alpha_{j}z^{j}\right)\left(\sum_{k=0}^{\infty}\beta_{k}w^{k}\right)=
\sum_{j=0}^{\infty}\sum_{k=0}^{\infty}\alpha_{j}\beta_{k}z^{j}w^{k}.
$$
In other words $\forall\,j,k\in\{0,1,2,3,\ldots\}$ we have $c_{jk}=\alpha_{j}\beta_{k}$ and condition 1 is equivalent
to: $\alpha_{j}\beta_{k}\ge 0$, $\forall\,j,k\in\{0,1,2,3,\ldots\}$. In particular, since $\alpha_{0},\,\beta_{0}\ne 0$
and since $\forall\,j\in\{0,1,2,3,\ldots\}$ we have $\alpha_{j}\beta_{0}\ge 0$ and $\forall\,k\in\{0,1,2,3,\ldots\}$
we have $\alpha_{0}\beta_{k}\ge 0$ it follows that we must have:
$$
\arg\alpha_{0}=\arg\alpha_{1}=\arg\alpha_{2}=\ldots=-\arg\beta_{0}=-\arg\beta_{1}=-\arg\beta_{2}=\ldots.
$$
We interpret $\arg 0$ as $\arg\alpha_{0}$ if $\alpha_{j_{0}}=0$ and as $\arg\beta_{0}$ if $\beta_{k_{0}}=0$
($\arg\alpha_{j_{0}}=\arg\alpha_{0}$, $\arg\beta_{k_{0}}=\arg\beta_{0}$). In particular, if $\forall\,j\in\{0,1,2,3,\ldots\}$,
$\alpha_{j}\ge 0$, then also $\forall\,k\in\{0,1,2,3,\ldots\}$, $\beta_{k}\ge 0$. In that case $G(b)=\Re G(b)$ for $b\ge 0$.
This, unfortunately, implies that conditions 1 and 2 contradict one another. Hence no such an entire function $G(w)$ exists.
Thus we have no factorization of $f(z)$ as $g(h(z))$ where $g^{'}(w)=e^{G(w)}$. The reason is the following: We have just seen
that condition 1 could be written as follows: $e^{-G(w)}=\sum_{k=0}^{\infty}\beta_{k}w^{k}$ where $\forall\,k\in\{0,1,2,3,\ldots\}$,
$\beta_{k}\ge 0$. Hence, if $G(w)$ is not a constant function (as must be the case in factorization), then
$\lim_{b\rightarrow\infty,\,b\in\mathbb{R}}e^{-G(b)}=\lim_{b\rightarrow\infty,\,b\in\mathbb{R}}\sum_{k=0}^{\infty}
\beta_{k}b^{k}=\infty$, for the coefficients $\beta_{k}\ge 0$ and at least two are non-zero. Thus: $\lim_{b\rightarrow\infty,\,
b\in\mathbb{R}}e^{G(b)}=0$, and since $G(w)$ is a non-constant entire function, we have, in fact, $0=\lim_{b\rightarrow\infty,\,
b\in\mathbb{R}}be^{G(b)}=\lim_{b\rightarrow\infty,\,b\in\mathbb{R}}be^{\Re G(b)}$. This violates condition 2. Thus such
a straight forward application of the majorant method of Cauchy, fails.
\end{remark}

\section{A recursion induced by factorization}
Let us assume that we have the following factorization in entire functions: $f(z)=g(h(z))$. We differentiate once
and obtain: $f^{'}(z)=g^{'}(h(z))h^{'}(z)$, and hence $f^{'}(z)/h^{'}(z)=g^{'}(h(z))$. This means that the quotient
$f^{'}(z)/h^{'}(z)$ is an entire function (and not merely a meromorphic function), and it shares the same right
factor $h(z)$ with $f(z)$. This motivates the following  definition of a recursive sequence:
\begin{equation}
\label{eq9}
f_{0}(z)=f(z),\,\,f_{n+1}(z)=\frac{f_{n}^{'}(z)}{h^{'}(z)}\,\,\,{\rm for}\,\,n\in\{0,1,2,3,\ldots\}.
\end{equation}
A simple inductive argument proves the following:

\begin{proposition}\label{prop9}
Let $f,\,g$ and $h$ be entire functions that satisfy the composition relation $f(z)=g(h(z))$. Let the sequence 
$\{f_{n}(z)\}_{n=1}^{\infty}$ be defined by the recursion in equation (\ref{eq9}). Then, \\
1. $\forall\,n\in\{0,1,2,3,\ldots\}$, $h^{'}(z)$ divides $f^{'}_{n}(z)$ within the entire functions. \\
2. $\forall\,n\in\{0,1,2,3,\ldots\}$, $f_{n}(z)=g^{(n)}(z)$. \\
Thus all the functions in the sequence $\{f_{0},f_{1},f_{2},f_{3},\ldots\}$ share the right factor $h(z)$, i.e.,
$f_{n}(z)\in R_{h}(E)=\{\phi\circ h\,|\,\phi\in E=\{{\rm entire}\,{\rm functions}\}\}$.
\end{proposition}

\begin{remark}\label{rem8}
1. Here are the first four elements in the sequence $\{f_{0},f_{1},f_{2},f_{3},\ldots\}$:
$$
f_{0}=f,\,f_{1}=\frac{f^{'}}{h^{'}},\,f_{2}=\frac{f^{''}h^{'}-f^{'}h^{''}}{h'^{3}},
$$
$$
f_{3}=\frac{f^{'''}h'^{2}-3f^{''}h^{'}h^{''}+f^{'}\cdot (3h^{''}-h^{'''}h^{'})}{h'^{5}}.
$$
2. In the case that $f\in {\rm elh}(\mathbb{C})$ and $f=g\circ h$, it follows (by the chain rule $f'=(g'\circ h)\cdot h'$)
that also $h\in {\rm elh}(\mathbb{C})$ and hence $h'(z)\ne 0$ $\forall\,z\in\mathbb{C}$, and we can divide by $h'(z)$. \\
3. We note that 
$$
f_{n+1}\cdot f_{n}=\frac{f'_{n}f_{n}}{h'}=\frac{\left(\frac{1}{2}f_{n}^{2}\right)'}{h'},
$$
where we have:
$$
\frac{1}{2}f(z)^{2}=\frac{1}{2}g(h(z))^{2}=g_{2}(h(z)),
$$
where $g_{2}=g^{2}/2$. So $f_{1}f_{0}$ also has a factorization $g_{2}(h(z))$ with the same right factor, $h(z)$.
\end{remark}

\section{A first order non-linear ordinary differential equation related to factorization over ${\rm elh}(\mathbb{C})$}

Our problem is the following: Given an entire function $F(z)$ such that $F(z)\not\equiv {\rm Const.}$, find two
entire functions $G(w)$ and $L(z)$ such that $G(w)\not\equiv {\rm Const.}$, $L(z)\not\equiv az+b$ for any $a,b\in\mathbb{C}$,
and such that:
\begin{equation}
\label{eq10}
\left\{\begin{array}{l} L{'}(z)=e^{F(z)}e^{-G(L(z))} \\ L(0)=0 \end{array}\right.,
\end{equation}
or prove that such entire functions $G(w)$ and $L(z)$ do not exist. This problem seems to be related to analytic continuation
(i.e. holomorphic extension). There are a variety of results on this topic. For example here is Theorem 17.0.6, on page 431
of \cite{im}. It is a classical result: \\
\\
{\bf Theorem.} (Painleve's Theorem) Let $E$ be a closed set of linear measure zero in a domain $\Omega\subseteq\mathbb{C}$.
Then every bounded holomorphic function $f\,:\,\Omega\rightarrow\mathbb{C}$ extends to a bounded holomorphic function
$f\,:\,\Omega\rightarrow\mathbb{C}$. In particular, if $\Omega=\mathbb{C}$, then $f$ is constant. \\
\\
By the standard theory of ode's we know the following: If we pick any entire function $G(w)$, $G(w)\not\equiv {\rm Const.}$,
then there exists a unique solution $L(z)$ for the system (\ref{eq10}), where $L(z)\in H(D(0,R(G)))$ for some $R(G)>0$, but
$\forall\,\epsilon>0$, $L(z)\not\in H(D(0,R(G)+\epsilon))$, unless $R(G)=+\infty$, i.e. unless $L(z)$ is an entire function,
in which case our problem is solved. Thus we will assume that $R(G)<\infty$, in which case we know that there is a point
$\alpha\in\partial D(0,R(G))$, beyond which the holomorphic function $L(z)$, can not be analytically continued.

\begin{theorem}\label{thm2}
Let $F(z)$ and $G(w)$ be entire functions such that $F(z),\,G(w)\not\equiv {\rm Const.}$ and let $L(z)$ be a solution
to the system (\ref{eq10}) such that $L(z)\in H(D(0,R(G)))$ for some $0<R(G)<\infty$, and $\forall\,\epsilon>0$,
$L(z)\not\in H(D(0,R(G)+\epsilon))$. Let $f(z)=\int_{0}^{z}e^{F(t)}dt$ and $g(w)=\int_{0}^{w}e^{G(w)}dw$ be the two
induced entire functions in ${\rm elh}(\mathbb{C})$. Then: \\
\\
(1) There is no sequence $\{z_{n}\}_{n=1}^{\infty}\subseteq D(0,R(G))$, such that $\lim z_{n}=\alpha\in\overline{D}(0,R(G))$,
$\{L(z_{n})\}_{n=1}^{\infty}$ is a bounded sequence but $\lim L(z_{n})$ does not exist. \\
\\
(2) There is a subset $A_{\infty}$ of $\partial D(0,R(G))$, such that $\forall\,\alpha\in\partial D(0,R(G))-A_{\infty}$ the
limit $\lim_{z\rightarrow\alpha,\,z\in D(0,R(G))}L(z)=a$ exists and is a complex number. \\
\\
(3) $\forall\,\alpha\in\partial D(0,R(G))-A_{\infty}$ there is an analytic continuation of $L(z)$ to a small enough disk,
$D(\alpha,\epsilon_{\alpha})$, $\epsilon_{\alpha}>0$. This, of course, implies part (2), but part (2) is needed in order
to prove this much stronger claim. \\
\\
(4) The set $A_{\infty}\subseteq\partial D(0,R(G))$ has the following properties:

(a) $A_{\infty}\ne\emptyset$.

(b) $\forall\,\alpha\in A_{\infty}$, $\lim_{z\rightarrow\alpha,\,z\in D(0,R(G))}|L(z)|=\infty$ 

and $\lim_{z\rightarrow\alpha,\,z\in D(0,R(G))} g(L(z))=f(\alpha)$.

(c) $f(A_{\infty})\subseteq A_{g}$ ($A_{g}$ is the set of all the finite asymptotic values of the 

entire function $g(w)$).

(d) $A_{\infty}$ is a closed subset of $\partial D(0,R(G))$.

(e) $|A_{\infty}|=0$, i.e. the one dimensional linear Lebesgue measure of $A_{\infty}$

is $0$.
\end{theorem}
\noindent
{\bf Proof.} \\
(1) Let us suppose that the claim is false and that $\{z_{n}\}_{n=1}^{\infty}$ is a counterexample.
Then $\lim z_{n}=\alpha$, $\{L(z_{n})\}_{n=1}^{\infty}$ is a bounded sequence, but $\lim L(z_{n})$ does not exist.
We will derive a contradiction. This will prove (1). By the assumptions on $F(z),\,G(w),\,L(z),\,f(z)$ and $g(w)$
we have the local factorization $f(z)=g(L(z))$, $\forall\,z\in D(0,R(G))$. We can find two sub-sequences
$\{a_{n}\}_{n=1}^{\infty}$ and $\{b_{n}\}_{n=1}^{\infty}$ of the convergent sequence $\{z_{n}\}_{n=1}^{\infty}$
such that $\lim L(a_{n})=a$, $\lim L(b_{n})=b$ and $a\ne b$. Let ${\rm Cl}_{L}(\alpha)$ be the set of all
the possible limits $\lim L(c_{n})$ over all of the sequences $\{c_{n}\}_{n=1}^{\infty}\subseteq D(0,R(G))$
such that $\lim c_{n}=\alpha$. Then $a,b\in {\rm Cl}_{L}(\alpha)$, $a\ne b$. The set ${\rm Cl}_{L}(\alpha)$
is called the cluster set of the function $L$ at $\alpha$. It is known that in our case ${\rm Cl}_{L}(\alpha)$
is a continuum. Hence $|{\rm Cl}_{L}(\alpha)|=\aleph_{1}$. For any sequence $\{c_{n}\}_{n=1}^{\infty}\subseteq
D(0,R(G))$, such that $\lim c_{n}=\alpha$ and $\lim L(c_{n})=c\in {\rm Cl}_{L}(\alpha)$ we have:
$$
f(c_{n})=g(L(c_{n}))\rightarrow g(c)=f(\alpha).
$$
Hence $g({\rm CL}_{L}(\alpha))=\{f(\alpha)\}$. Since $g(w)$ is an entire function, we deduce that $g(w)\equiv f(\alpha)$.
But $g(w)\in {\rm elh}(\mathbb{C})$ and in particular can not be constant. This is a contradiction. \\
\\
(2) Let us denote by $A_{\infty}$ the set complementary with respect to $\partial D(0,R(G))$ of the set
of points $\alpha$ we dealt with in part (1). By what we proved in part (1), $\forall\,\alpha\in\partial D(0,R(G))-A_{\infty}$
we have no two sequences $\{a_{n}\}_{n=1}^{\infty}$ and $\{b_{n}\}_{n=1}^{\infty}$ contained in $D(0,R(G))$ so that
$\lim a_{n}=\lim b_{n}=\alpha$, but the two limits $a=\lim L(a_{n})$ and $b=\lim L(b_{n})$ exist but differ from
one another, $a\ne b$. Thus, indeed, for any such an $\alpha$ we have $\lim_{z\rightarrow,\,z\in D(0,R(G))}L(z)=a$,
the same limit. We only need to remark that this is also true when one limit, say $a$ is a finite complex number
but $|b|=+\infty$. The argument in the proof of part (1) can be adjusted to imply also this case. \\
\\
(3) Let us take an $\alpha\in\partial D(0,R(G))-A_{\infty}$. Then by part (2) we have 
$\lim_{z\rightarrow\alpha,\,z\in D(0,R(G))} L(z)=a\in\mathbb{C}$. Let us look at the following system:
\begin{equation}
\label{eq11}
\left\{\begin{array}{l} L^{'}(z)=e^{F(z)}e^{-G(L(z))} \\ L(\alpha)=a \end{array}\right..
\end{equation}
The general theory of ode's imply that there is a positive radius $\epsilon_{\alpha}>0$, such that in the disk
$D(\alpha,\epsilon_{\alpha})$ the system (\ref{eq11}) has a unique holomorphic solution $L(z)$. We used the same
notation $L(z)$ as for the solution of the system (\ref{eq10}) because this is indeed an analytic continuation
of the original $L(z)$ to the disk $D(\alpha,\epsilon_{\alpha})$. \\
\\
(4) (a) If $A_{\infty}=\emptyset$, then by part (3) it follows that the solution $L(z)$ of the system (\ref{eq10})
belongs to $H(D(0,R(G)+\epsilon))$ for some $\epsilon>0$. This contradicts our assumptions on $L(z)$ and on
$R(G)<\infty$. \\
\\
(b) If $\alpha\in A_{\infty}$, then by parts (1) and (2) it follows that for any sequence $\{a_{n}\}_{n=1}^{\infty}$
in $D(0,R(G))$ such that $\lim a_{n}=\alpha$, we must have $\lim |L(a_{n})|=\infty$. We have a second interpretation
of our functions, not as related by an ode, but by composition of mappings: $f(z)=g(L(z))$. So: $f(a_{n})=g(L(a_{n}))
\rightarrow f(\alpha$. Hence $\lim_{z\rightarrow\alpha,\,z\in D(0,R(G))}g(L(z))=f(\alpha)$. \\
\\
(c) This follows by the definition of the asymptotic set, $A_{g}$, of the entire function $g(w)$, and by part (4)(b). \\
\\
(d) $A_{\infty}$ is closed because its complementary set with respect to $\partial D(0,R(G))$ is an open subset of
$\partial D(0,R(G))$ (by part (3)). \\
\\
(e) By part (4)(b) we have: $\forall\,\alpha\in A_{\infty}$, $\lim_{z\rightarrow\alpha,\,z\in D(0,R(G))}|L(z)|=\infty$.
By standard results on the boundary values of holomorphic functions we deduce that if $|A_{\infty}|>0$, then $L(z)|\equiv\infty$,
$\forall\,z\in D(0,R(G))$. $\qed $ \\
\\
Theorem \ref{thm2} deals with non-entire solution $L(z)$ of the system (\ref{eq10}). It gives the boundary values
of the solution on the maximal disk in which $L(z)$ is holomorphic. The disk is centered at the point at which the initial
value condition is given. In the theorem we used the notation $D(0,R(G))$. The assumption was that the functions $F(z)$
and $G(w)$ are fixed entire functions. However, these two mappings have non identical roles. Namely the entire function
$F(z)$ is given in the system (\ref{eq10}). On the other hand the functions $G(w)$ and $L(z)$ are entire functions that
satisfy $G(w)\not\equiv {\rm Const.}$ and $L(z)\not\equiv az+b$ for any $a,b\in\mathbb{C}$. Our method is to fix (arbitrarily)
an entire function $G(w)$ (subject only to the restriction $G(w)\not\equiv {\rm Const.}$) and see how large could the
maximal disc of holomorphy of $L(z)$ be. It's radius clearly depends on $G(w)$ and hence the notation $R(G)$. We also
investigated the obstacles that prohibit $L(z)$ from being an entire function, i.e. the nature the singularities of $L(z)$
on the boundary $\partial D(0,R(G))$. This set was denoted in Theorem \ref{thm2} by $A_{\infty}$. Also $A_{\infty}$
depends on $G(w)$, and we could have denoted it by $A_{\infty}(G)$. This set of singularities is described in parts
(4)(a)(b)(c)(d). It is characterized by the following two limits:
$$
\alpha\in A_{\infty}\Leftrightarrow \lim_{z\rightarrow\alpha,\,z\in D(0,R(G))}|L(z)|=\infty\,\,{\rm and}\,\,
\lim_{z\rightarrow\alpha,\,z\in D(0,R(G))}g(L(z))=f(\alpha),
$$
where $f(z)=\int_{0}^{z}e^{F(t)}dt$ and $g(w)=\int_{0}^{w}e^{G(t)}dt$. Thus $f(A_{\infty})\subseteq A_{g}$. We also
proved that $A_{\infty}$ is a closed subset of $\partial D(0,R(G))$ with a one dimensional linear Lebesgue
measure $0$, i.e. $|A_{\infty}|=0$. Our next intention is to try and understand the way in which the radius $R(G)$ or
the set of singularities on $L(z)$, $A_{\infty}(G)$ depend on $G(w)$. Hopefully there exists an appropriate
perturbation $G\rightarrow G_{0}$ such that $r(G)\rightarrow R(G_{0})=\infty$, or, equivalently $A_{\infty}(G)\rightarrow
A_{\infty}(G_{0})=\emptyset$. For that purpose we define now
a maximal analytic continuation of $L(z)$ (for a fixed $G(w)$). A main ingredient will be the trick introduced in the proof
of part (3), that of "moving" the system (\ref{eq10}) to a new location centered at a point $\alpha\in\partial D(0,R(G))-A_{\infty}$.
This perturbation of the system (\ref{eq10}) changes only the center and the initial value. The center, as mentioned is $\alpha$, 
and the initial value is determined by the limit $a=\lim_{z\rightarrow\alpha,\,z\in D(0,R(G))}L(z)\in\mathbb{C}$. We will
now take the largest possible radius of holomorphy of the solution $L(z)$ of the perturbed system (\ref{eq11}). This new largest
disk of holomorphy can be denoted by $D(\alpha_{1},R_{1}(G))$. We will call it a generation $1$ disk. Any point
of $\partial D(0,R(G))-A_{\infty}$ is the center of a generation $1$ disk. Along these lines we will call the first
maximal disk $D(0,R(G))$, the generation $0$ disk and write sometimes $D(0,R(G))=D(0,R_{0}(G))$. The singular set of the
solution $L(z)$ of a generation $1$ system (\ref{eq11}), is a closed subset of $\partial D(\alpha_{1},R_{1}(G))$, which
is non-empty and of one dimensional linear Lebesgue measure $0$. This generation $1$ singular set can be denoted by
$A_{1,\infty}$. The generation $0$ singular set is the non-empty set $A_{\infty}=A_{0,\infty}$ which is a subset of
$\partial D(0,R_{0}(G))$ and which was described in Theorem \ref{thm2}. Note that also the generation $1$ singular set
is characterized by two limits:
$$
\alpha\in A_{1,\infty}\Leftrightarrow\lim_{z\rightarrow\alpha,\,z\in D(0,R_{1}(G))}|L(z)|=\infty\,\,{\rm and}\,\,
\lim_{z\rightarrow\alpha,\,z\in D(0,R_{1}(G))}g(L(z))=f(\alpha).
$$
Also $A_{1,\infty}$ is a closed subset of $\partial D(\alpha_{1},R_{1}(G))$ which has one dimensional linear Lebesgue
measure $0$. The process can now be continued indefinitely (but only $\aleph_{0}$ times because of area considerations)
to produce any generation of analytic continuation in the the maximal disk of holomorphy $D(\alpha_{k},R_{k}(G))$ where
$\alpha_{k}\in \partial D(0,R_{k-1}(G))-A_{k-1,\infty}$. If we denote: $a_{k}=\lim_{z\rightarrow\alpha_{k},\,z\in
D(\alpha_{k-1},R_{k-1}(G))}L(z)$, where by an application of Theorem \ref{thm2} part (2) the last limit exists and is
a finite complex number, then the $k$'th generation initial value problem is:
\begin{equation}
\label{eq12}
\left\{\begin{array}{l} L^{'}(z)=e^{F(z)}e^{-G(L(z))} \\ L(\alpha_{k})=a_{k} \end{array}\right.,\,\,k=0,1,2,3,\ldots.
\end{equation}
It defines $L(z)$ uniquely in $D(\alpha_{k},R_{k}(G))$ so that $0<R_{k}(G)<\infty$, and $\forall\,\epsilon>0$,
$L(z)\not\in D(\alpha_{k},R_{k}(G)+\epsilon)$. The totality of the singular sets of any generation, namely 
$$
{\rm sing}(L)=\overline{\bigcup_{k=0}^{\infty}\left(\cup A_{k,\infty}\right)},
$$
is the singular set of $L(z)$. We note that ${\rm sing}(L)$ is a non-empty closed subset of the maximal Riemann surface
that uniformizes $L(z)$. It is not clear if it has linear Lebesgue measure $0$. It is the closure of countably many sets 
of measure $0$ and hence in principle might have non zero measure. The point is that the geometric structure of those sets 
if very special. They lie on the boundaries of the maximal holomorphy disks, and the contribution of the $k$'th generation,
i.e. the sets $A_{k,\infty}$ is outside the the union of all the maximal holomorphy disks up to and including those of 
generation $k-1$. Also, any point:
$$
\alpha\in\bigcup_{k=0}^{\infty}\left(\cup A_{k,\infty}\right),
$$
(not the closure of this set!) lies on the boundary of some finite generation maximal disk of holomorphy, $\alpha\in
\partial D(\alpha_{k},R_{k}(G))$. In fact, $\alpha\in A_{k,\infty}$ and so is characterized by the two limits characterization:
$$
\lim_{z\rightarrow\alpha,\,z\in D(\alpha_{k},R_{k}(G))}|L(z)|=\infty\,\,{\rm and}\,\,
\lim_{z\rightarrow\alpha,\,z\in D(\alpha_{k},R_{k}(G))}g(L(z))=f(\alpha).
$$
Thus $f({\rm sing}(L))\subseteq A_{g}$ (where the asymptotic tract is taken on one of the branches of the total Riemann surface
that uniformizes $L(z)$).

\begin{theorem}\label{thm3}
The set ${\rm sing}(L)$ of all the singular points of a maximal holomorphic extension, $L(z)$ of the solution $L(z)$ of the 
system (\ref{eq10}), is a perfect subset of the uniformizing Riemann surface.
\end{theorem}
\noindent
{\bf Proof.} \\
Suppose that $\alpha\in {\rm sing}(L)$ is an isolated point (of ${\rm sing}(L)$). Then since ${\rm sing}(L)$ is the closure
(on the uniformizing Riemann surface of $L(z)$) of the countable union $\bigcup_{k=0}^{\infty}\left(\cup A_{k,\infty}\right)$
where each $A_{k,\infty}$ is a closed and null subset of a circle $\partial D(\alpha_{k},R_{k}(G))$, it follows that
$\alpha\in\bigcup_{k=0}^{\infty}\left(\cup A_{k,\infty}\right)$ and does not belong to the portion of the derived set
of $\bigcup_{k=0}^{\infty}\left(\cup A_{k,\infty}\right)$ which lies outside the set 
$\bigcup_{k=0}^{\infty}\left(\cup A_{k,\infty}\right)$ itself. Hence $\exists\,k\in\mathbb{Z}^{+}\cup\{0\}$ such that
$\alpha\in A_{k,\infty}$. We conclude that $\alpha$ is an isolated singular point of the holomorphic mapping $L(z)$, and
also that $\lim_{z\rightarrow\alpha,\,z\in D(\alpha_{k},R_{k}(G))}|L(z)|=\infty$. By the classification theorem of 
isolated singularities of an holomorphic function, we deduce that either $\alpha$ is a pole of $L(z)$ of some
order $k_{0}$, or $\alpha$ is a principle singularity of $L(z)$. If $\alpha$ is such a pole of $L(z)$, then there
is a punctured neighborhood of $\alpha$, $V(\alpha)$ and a bounded holomorphic function $m(z)\in H(V(\alpha)\cup\{\alpha\})$,
such that $m(\alpha)\ne 0$, so that $\forall\,z\in V(\alpha)$, $L(z)=m(z)\cdot (z-\alpha)^{-k_{0}}$. Here 
$k_{0}\in\mathbb{Z}^{+}$. We recall that the entire function $g(w)=\int_{0}^{w}e^{G(u)}du$ is in ${\rm elh}(\mathbb{C})$
and can not be a polynomial, so $g(w)$ is a a transcendental entire function. Hence by $f(z)=g(L(z))=
g(m(z)\cdot (z-\alpha)^{-k_{0}})$, $\forall\,z\in V(\alpha)$ it follows that $f(z)$ has an essential singularity
at $z=\alpha$. This contradicts the fact that $f(z)$ is an entire function. If, on the other hand, $z=\alpha$ is 
an essential singular point of $L(z)$, then the composition $g(L(z))$ has an essential singular point at $z=\alpha$
(we recall that $g(w)\in {\rm elh}(\mathbb{C})$), but $f(z)=g(L(z))$ so we obtain the same contradiction. Now
the assertion is proved. $\qed $

\begin{theorem}\label{thm4}
If $f\in {\rm elh}(\mathbb{C})$ has a Picard exceptional value, then $f$ is not a prime entire function.
\end{theorem}
\noindent
{\bf Proof.} \\
Let us suppose that the number $a\in\mathbb{C}$ is a Picard exceptional value of $f(z)$. This means that
$f(\mathbb{C})=\Im\{ f\}=\mathbb{C}-\{a\}$. Equivalently, this means that the function $f(z)-a$ is an entire
function that has no zeroes. Since $\mathbb{C}$ is a simply connected domain, it follows that we can define
$L(z)=\log(f(z)-a)$ as an entire function. By the assumption that $f$ has a Picard exceptional value, it
follows that $f(z)\not\equiv c\cdot z+d$ for any $c\in\mathbb{C}^{\times}$ and $d\in\mathbb{C}$. Thus $f(z)$
is a transcendental entire function. If the entire function $L(z)$ is an affine function, say $L(z)=c\cdot z+d$,
then $c\cdot z+d=\log(f(z)-a)$. So $f(z)=a+e^{d}e^{c\cdot z}$ and this function is not a prime entire function.
For example, we have $e^{z}=z^{2}\circ e^{z/2}$, $e^{c\cdot z}=z^{2}\circ e^{c\cdot z/2}$ and hence:
$$
f(z)=(a+e^{d}z)\circ (z^{2}\circ e^{c\cdot z/2}),
$$
a non trivial factorization of $f(z)$. On the other hand, if $L(z)$ is not an affine function, then it is
a transcendental entire function. The reason is that $L^{'}(z)=f^{'}(z)/(f(z)-a)$ and so $L^{'}(z)\ne 0$
$\forall\,z\in\mathbb{C}$. Thus $L(z)\in {\rm elh}(\mathbb{C})- {\rm Aut}(\mathbb{C})$ which implies that
$L(z)$ is a transcendental function. By $f(z)=a+e^{L(z)}$ it follows that $f(z)=(a+z)\circ e^{z}\circ L(z)$,
a non trivial factorization of $f(z)$. $\qed $ \\
\\
Thus if $f(z)\in {\rm elh}(\mathbb{C})$ is a prime entire function, then $f(z)$ has no Picard exceptional 
value. This means that $f(\mathbb{C})=\mathbb{C}$. Using the result on page 251 of the book \cite{ru} we
know that a Picard exceptional value of an entire function (i.e. an omitted value of an entire function)
is an asymptotic value. Thus in our case, when $f\in {\rm elh}(\mathbb{C})$, $f(\mathbb{C})=\mathbb{C}$, has
asymptotic values, then these are not omitted values of $f$.

In order to find out if ${\rm elh}(\mathbb{C})$ contains primes, it follows, by Theorem \ref{thm4} that it is
sufficient to consider only functions $f(z)\in {\rm elh}(\mathbb{C})$, such that $f(\mathbb{C})=\mathbb{C}$.
For such a function, if there is a factorization $f(z)=g(L(z))$ by entire functions then necessarily $g(z)$
has no omitted values, i.e. $g(\mathbb{C})=\mathbb{C}$. Also by $f^{'}(z)=L^{'}(z)g^{'}(L(z))$, it follows
that $L^{'}(z)\ne 0$ $\forall\,z\in\mathbb{C}$ and that $g^{'}(z)\ne 0$ $\forall\,z\in L(\mathbb{C})$. Thus
$L\in {\rm elh}(\mathbb{C})$. Now, either $L(\mathbb{C})=\mathbb{C}$ in which case also $g(z)\in {\rm elh}(\mathbb{C})$,
or $L(\mathbb{C})=\mathbb{C}-\{a\}$, in which case $g(z)=\int_{0}^{z}(t-a)^{N}e^{G(t)}dt$ for some
$N\in \mathbb{Z}^{+}\cup\{0\}$ and an entire function $G(t)$. Using the equations before equation (\ref{eq3}), we
note that when $L(\mathbb{C})=\mathbb{C}-\{a\}$, then
$$
e^{F(z)}=e^{H(z)}\cdot\left(\int_{0}^{z}e^{H(w)}dw-a\right)^{N}\cdot e^{G(\int_{0}^{z}e^{H(w)}dw)},
$$
where we have: $f(z)=\int_{0}^{z}e^{F(w)}dw$, $g(z)=\int_{0}^{z}(w-a)^{N}e^{G(w)}dw$ and $L(z)=\int_{0}^{z}e^{H(w)}dw$.
So instead of the system (\ref{eq10}) or the system (\ref{eq11}) we have the following system:
$$
\left\{\begin{array}{l} L^{'}(z)=e^{F(z)}(L(z)-a)^{-N}e^{-G(L(z))} \\ L(\alpha)=a_{0}\end{array}\right..
$$
Next we recall that when $f(z)=g(h(z))$ for entire functions, then there are some quantitative restrictions
among the orders of growth of these functions, or among their radial maximum modulus functions. We will quote
results du to Clunie and to P\'olya. We refer the reader to the book \cite{cc} (See also \cite{ha}). On page 207 we find the
following: \\
\\
{\bf Theorem A.2.} (Clunie) Let $f(z)$ and $g(z)$ be entire functions with $g(0)=0$. Let $\rho$ satisfy
$0<\rho<1$ and let $c(\rho)=\left(\frac{1-\rho^{2}}{4\rho}\right)$. Then for $R\ge 0$,
$$
M(R,f\circ g)\ge M\left(c(\rho)M(\rho R,g),f\right).
$$
On Page 208 we find: \\
\\
{\bf Corollary A.1.} (P\'olya) Let $f(z)$, $g(z)$ and $h(z)$ be entire functions with $h(z)=f(g(z))$. If $g(0)=0$,
then there exists an absolute constant $c$, $0<c<1$ such that for all $r>0$ the following inequalities holds:
$$
M(r,h)\ge M\left(cM\left(\frac{r}{2},g\right),f\right).
$$
On page 209 we find: \\
\\
{\bf Theorem A.3.} If $f(z)$ and $g(z)$ are two entire functions such that $f(g)$ is of finite order (lower order),
then \\
(i) either $g(z)$ is a polynomial and $f(z)$ is of finite order (lower order), \\
or \\
(ii) $g(z)$ is not a polynomial but a function of finite order (lower order) and $f(z)$ is of zero order (lower order). \\
\\
On page 210: \\
\\
{\bf Theorem A.4.} Let $f$ and $g$ be two transcendental entire functions. Then:
$$
\lim_{r\rightarrow\infty}\frac{\log M(r,f\circ g)}{\log M(r,f)}=\infty,
$$
and
$$
\lim_{r\rightarrow\infty}\frac{\log M(r,f\circ g)}{\log M(r,g)}=\infty.
$$
On page 213: \\
\\
{\bf Theorem A.7} (Edrei and Fuchs) Let $f(z)$ be an entire function that is not of zero order and $g(z)$ be
a transcendental entire function, then $f(g)$ is of infinite order. \\
\\
We recall that by Theorem \ref{thm2} part 4(c) it follows that $f(A_{\infty})\subseteq A_{g}$. In particular, if it
were possible to find a function $g\in {\rm elh}(\mathbb{C})-{\rm Aut}(\mathbb{C})$, for which $A_{g}$ is small
(for example $A_{g}=\emptyset$), then necessarily the singular locus, $A_{\infty}$ of $L(z)$ on $\partial D(0,R(G))$
had to be small (say $A_{\infty}=\emptyset$). We know that the asymptotic variety $A_{g}$ tends to be small when the
order of $g$ is small. For example, by the well known result of the Wiman-Valiron theory, if the order $\rho(u)$ of an 
entire function $u(z)$ satisfies $\rho(u)<1/2$, then $A_{u}=\emptyset$. However, the well known result of J. Hadamard,
\cite{h}, implies that for any $g\in {\rm elh}(\mathbb{C})-{\rm Aut}(\mathbb{C})$ we must have $A_{g}\ne\emptyset$. Thus
by the result mentioned above of the Wiman-Valiron theory it follows immediately that if $g\in {\rm elh}(\mathbb{C})-
{\rm Aut}(\mathbb{C})$ then necessarily $\rho(g)\ge 1/2$. In fact, one can state more accurate estimates along these
lines. We refer to the paper \cite{e}. On page 8 we find the following result of Clunie, Eremenko, Langley and Rossi: \\
\\
{\bf Theorem 3. (\cite{e})} (Clunie, Eremenko, Langley and Rossi) Let $f$ be a transcendental meromorphic function
of order $\rho$. \\
(a) If $\rho<1$, then $f^{'}$ has infinitely many zeroes. \\
(b) If $\rho<1/2$, then $f^{'}/f$ has infinitely many zeroes. \\
(c) If $f$ is entire, and $\rho<1$, then $f^{'}/f$ has infinitely many zeroes. \\
\\
So the order of growth $\rho$ of functions in ${\rm elh}(\mathbb{C})-{\rm Aut}(\mathbb{C})$ can not be smaller than $1$.
We now refer to the paper \cite{sm} for results on asymptotic values of functions in ${\rm elh}(\mathbb{C})-{\rm Aut}(\mathbb{C})$.
Firstly, few examples for functions $g\in {\rm elh}(\mathbb{C})-{\rm Aut}(\mathbb{C})$ such that $g(\mathbb{C})=\mathbb{C}$.
We refer to the two examples in that paper, on page 640: 

\begin{example}\label{exmp1}
Set:
$$
h(z)=\sum_{n=1}^{\infty}\frac{z^{n}}{n\cdot n!}=\int_{0}^{z}\frac{e^{\xi}-1}{\xi}d\xi.
$$
Then $h$ is a transcendental entire function and $1+zh^{'}(z)=e^{z}$. The function $g(z)=z\cdot e^{h(z)}$ is also
a transcendental entire function and $g^{'}(z)=\exp(z+h(z))$. In particular, $g\in {\rm elh}(\mathbb{C})-{\rm Aut}(\mathbb{C})$.
Because $g$ has an essential singularity at $\infty$, Picard's Big Theorem assures that $g$ assumes every finite complex
value infinitely often with at most one exception. Since $g$ takes the value $0$ exactly once, every nonzero complex
number is assumed infinitely often by $g$. Therefore $g(\mathbb{C})=\mathbb{C}$.
\end{example}

\begin{example}\label{exmp2} (Attributed to E. Calabi) Let $g(z)=\int_{0}^{z}e^{-\xi^{2}}d\xi$, then $g^{'}(z)=e^{-z^{2}}$, so
$g\in {\rm elh}(\mathbb{C})-{\rm Aut}(\mathbb{C})$, $g(0)=0$ and $g$ is an odd function. If $g$ omits the value $w_{0}$,
then $w_{0}\ne 0$ and by Picard's Theorem $g$ cannot omit $-w_{0}$. But $g(z_{0})=-w_{0}$ implies that $g(-z_{0})=w_{0}$,
a contradiction. Thus $g(\mathbb{C})=\mathbb{C}$.
\end{example}
\noindent
The following is a natural question:

\begin{problem}\label{prb1}
Is any of the functions:
$$
z\cdot\exp\left(\int_{0}^{z}\frac{e^{\xi}-1}{\xi}d\xi\right),\,\,\,\int_{0}^{z}e^{-\xi^{2k}}d\xi,\,\,{\rm for}\,k\in\mathbb{Z}^{+},
$$
a prime entire mapping?
\end{problem}
\noindent
Next, we recall three results that are given on page 641 of \cite{sm}. \\
\\
{\bf Proposition. (\cite{sm})} Suppose that $g\in {\rm elh}(\mathbb{C})$. Given $a\in\mathbb{C}$ set $b=g(a)$. Let
$f_{a}$ denote the branch of $g^{-1}$ which is defined in a neighborhood of $b$ and satisfies $f_{a}(b)=a$. The
radius of convergence of the Taylor series expansion of $f_{a}$ about $b$ is designated by $r_{a}$. \\
(i) If $r_{a}=\infty$, then $g\in {\rm Aut}(\mathbb{C})$. \\
(ii) If $r_{a}<\infty$, then every singular point on the circle of convergence $\{ w\,|\,|w-b|=r_{a}\}$ of $f_{a}$
is an asymptotic value of $g$. \\
\\
{\bf Corollary 1. (\cite{sm})} If $g\in {\rm elh}(\mathbb{C})$, then either $g\in {\rm Aut}(\mathbb{C})$ or $g$ has
at least one finite asymptotic value (i.e. $A_{g}\ne\emptyset$). \\
\\
{\bf Corollary 2. (\cite{sm})} If $g\in {\rm elh}(\mathbb{C})$ and $g(\mathbb{C})=\mathbb{C}$, then either
$g\in {\rm Aut}(\mathbb{C})$ or $g$ has at least two finite asymptotic values (i.e. $|A_{g}|\ge 2$). \\
\\
Finally, on page 642 of \cite{sm} we find: \\
\\
{\bf Proposition. (\cite{sm})} (i) An entire function $g\in {\rm elh}(\mathbb{C})$ belongs to ${\rm Aut}(\mathbb{C})$
if and only if $g$ has no finite asymptotic value (i.e. $A_{g}=\emptyset$). \\
(ii) $g\in {\rm elh}(\mathbb{C})$ and $A_{g}=\{\alpha\}$ ($\alpha\in\mathbb{C}$) if and only if there are
$a\in\mathbb{C}^{\times}$ and $b\in\mathbb{C}$  such that $g(z)=\alpha+\exp(a\cdot z+b)$. \\
\\
\begin{remark}\label{rem9}
Part (i) of the last proposition follows by the theorem of J. Hadamard, \cite{h}, on local diffeomorphisms
$F\,:\,\mathbb{R}^{n}\rightarrow\mathbb{R}^{n}$.
\end{remark}

\begin{example}\label{exmp3}
Next, we refer to page 268 of the book, \cite{n}. Let $f(z)=\int_{0}^{z}e^{e^{t}}dt$. It's Picard exceptional value
$\infty$ is of deficiency $\delta(\infty)=1$. It has infinitely many finite asymptotic values, all of which have
deficiency zero, so the total deficiency reduces to $1$, and does not make it to the maximal value $2$. The reduction
in the total deficiency in this example is a consequence of the fact that the finite asymptotic values are all approximated
by the function $f(z)$ at asymptotically the same rate. Because of this symmetry, the asymptotic values (the finite ones)
must all have equal deficiencies. Since their number is infinite, while the total deficiency is bounded by $2$, these
deficiencies must all vanish. In particular there is no finite Picard exceptional value of $f(z)$. In other words
$f(\mathbb{C})=\mathbb{C}$.
\end{example}
\noindent
We recall that we would like to know if ${\rm elh}(\mathbb{C})$ contains prime entire mappings. By Theorem \ref{thm4}
it follows that such mappings have no Picard exceptional value. At this point using Example \ref{exmp2} of E. Calabi
we are close to answer that question.

\begin{theorem}\label{thm5}
The function $f(z)=\int_{0}^{z}e^{-t^{2}}dt$ belongs to ${\rm elh}(\mathbb{C})-{\rm Aut}(\mathbb{C})$. It satisfies
$f(\mathbb{C})=\mathbb{C}$. It has no non degenerate factorization of the type $f(z)=g(L(z))$, where
$L(z)\in {\rm elh}(\mathbb{C})-{\rm Aut}(\mathbb{C})$ satisfies $L(\mathbb{C})=\mathbb{C}$, and $g(z)\in {\rm elh}(\mathbb{C})
-{\rm Aut}(\mathbb{C})$.
\end{theorem}

\begin{remark}\label{rem10}
Clearly $g(z)$ must also satisfy $g(\mathbb{C})=\mathbb{C}$.
\end{remark}
\noindent
{\bf Proof.} \\
The function $f(z)$ is the construction of E. Calabi we described in Example \ref{exmp2}. Thus the only
thing we need to prove is the claim that $f(z)$ has no non degenerate factorization of the type above.
In the case that the inner factor $L(z)$ satisfies $L(\mathbb{C})=\mathbb{C}$, the system which is equivalent
to this factorization is of the type in equation (\ref{eq10}) or more generally in equation (\ref{eq11}), i.e.:
$$
\left\{\begin{array}{l} L^{'}(z)=e^{F(z)}e^{-G(L(z))} \\ L(\alpha)=a \end{array}\right..
$$
Thus the differential relation is $L^{'}(z)=e^{F(z)}e^{-G(L(z))}$, where the entire function $G(w)\not\equiv {\rm Const.}$.
We note that by $L\in {\rm elh}(\mathbb{C})$, it follows that $\log L^{'}(z)$ is an entire function and the differential
equation above can also be written as follows:
$$
\log L^{'}(z)=F(z)-G(L(z)).
$$
If $L(z)$ is an entire function of a finite order, then also $L^{'}(z)$ is an entire function and of the same order
(as that of $L(z)$). This can be found in Theorem 1.7 on page 115 of the book \cite{s}. A stronger theorem, is Theorem 2.1
on page 120 of \cite{s}. Since the order of growth of the composition $G(L(z))$ is at least as large as that of
$L(z)$ ($G$ is a non constant entire function), and since the order of $\log L^{'}(z)$ equals that of $\log L(z)$
which is strictly smaller than that of $L(z)$, the difference of the entire functions $F(z)-G(L(z))$ produces
cancellation in the order of growth of $G(L(z))$, with an error term which equals $\log L^{'}(z)\approx\log L(z)$. Coming
back to the function that was suggested by E. Calabi, we have $F(z)=-z^{2}$. Thus
$$
-z^{2}-G(L(z))=\log L^{'}(z)\approx\log L(z).
$$
We remark that this implies that $G(L(z))$ has zero order. By Theorem A.3 (quoted above from \cite{cc}) it follows
that: \\
(i) either $L(z)$ is a polynomial and $G(z)$ is of zero order, \\
or \\
(ii) $L(z)$ is not a polynomial but a function of zero order and also $G(z)$ is of zero order. \\
Since the only polynomials in ${\rm elh}(\mathbb{C})$ are the functions in ${\rm Aut}(\mathbb{C})$, but $L(z)$ is not
in ${\rm Aut}(\mathbb{C})$, it follows that in our situation only case (ii) is possible. By Theorem A.4 in \cite{cc}
it follows that if $G(z)$ and $L(z)$ are two transcendental entire functions, then:
$$
\lim_{r\rightarrow\infty}\frac{\log M(r,G\circ L)}{\log M(r,G)}=\lim_{r\rightarrow\infty}\frac{\log M(r,G\circ L)}
{\log M(r,L)}=\infty.
$$
In our case (assuming that $L(z)\in {\rm elh}(\mathbb{C})-{\rm Aut}(\mathbb{C})$) we have the following equations:
$$
\lim_{r\rightarrow\infty}\frac{M(r,G(L(z)))}{r^{2}}=\lim_{r\rightarrow\infty}\frac{M(r,G(L(z)))}{M(r,\log L^{'}(z)}=\infty.
$$
Note these equations do not include the logarithms. Hence, the equation $-z^{2}-G(L(z))=\log L^{'}(z)$ can not have
a solution where $G(w)\not\equiv {\rm Const.}$ is an entire function and where $L(z)\in {\rm elh}(\mathbb{C})-{\rm Aut}(\mathbb{C})$.
This concludes the proof of the theorem $\qed $ \\
\\
It is clear that the equation:
$$
lim_{r\rightarrow\infty}\frac{M(r,G(L(z)))}{r^{2}}=\infty,
$$
can be generalized as follows:
$$
\lim_{r\rightarrow\infty}\frac{M(r,G(L(z)))}{r^{N}}=\infty,\,\,\,\forall\,N\in\mathbb{Z}^{+}.
$$
Thus, in fact we have proved the following more general:

\begin{theorem}\label{thm6}
The function $f(z)=\int_{0}^{z} e^{p(t)}dt$, where $p(t)\in\mathbb{C}[t]-\mathbb{C}$, belongs to 
${\rm elh}(\mathbb{C})-{\rm Aut}(\mathbb{C})$. It has no non-degenerate factorization of the type
$f(z)=g(L(z))$, where $L(z)\in {\rm elh}(\mathbb{C})-{\rm Aut}(\mathbb{C})$ satisfies $L(\mathbb{C})=\mathbb{C}$,
and where $g(z)\in {\rm elh}(\mathbb{C})-{\rm Aut}(\mathbb{C})$.
\end{theorem}

\begin{remark}\label{rem11}
The claim in Theorem \ref{thm6} can not be strengthened to the following: $f(z)$ is a prime entire function. The
exponential function shows that. I.e. if $p(t)\equiv t$ then $f(z)=\int_{0}^{z}e^{t}dt=e^{z}-1=(z^{2}-1)\circ e^{z/2}$.
We note that in this example $L(z)=e^{z/2}$ has a Picard exceptional value.
\end{remark}
\noindent
Clearly, we now investigate the last possible case, namely that of which we allow the inner factor 
$L(z)\in {\rm elh}(\mathbb{C})-{\rm Aut}(\mathbb{C})$ to have a Picard exceptional value, say $a$.

\begin{theorem}\label{thm7}
The function $f(z)=\int_{0}^{z}e^{p(t)}dt$, where $p(t)\in\mathbb{C}[t]-\mathbb{C}$, belongs to
${\rm elh}(\mathbb{C})-{\rm Aut}(\mathbb{C})$. It is an entire prime mapping if and only if $f(\mathbb{C})=\mathbb{C}$.
\end{theorem}
\noindent
{\bf Proof.} \\
By Theorem \ref{thm6} we know that $f(z)$ has no non degenerate factorization of the type $f(z)=g(L(z))$ where $L(z)$
and $g(z)$ belong to ${\rm elh}(\mathbb{C})-{\rm Aut}(\mathbb{C})$, and where $L(z)$ has no Picard exceptional value,
i.e. $L(\mathbb{C})=\mathbb{C}$. It remains to deal with the possibility that $f(z)$ has a non degenerate factorization
in which $L(\mathbb{C})=\mathbb{C}-\{ a\}$. In this case the differential relation that characterizes this situation is:
$$
e^{p(z)}=L^{'}(z)\cdot(L(z)-a)^{N}\cdot e^{G(L(z))},
$$
where as always:
$$
f(z)=\int_{0}^{z}e^{p(t)}dt+c,\,\,\,g(z)=\int_{0}^{z}e^{G(t)}\cdot (t-a)^{N}dt+d\,\,\,\,\,{\rm where}\,\,N\in
\mathbb{Z}^{+}\cup\{ 0\},
$$
and where $G$ is entire which is non constant if $N=0$, but which might be constant if $N\in\mathbb{Z}^{+}$. Thus:
$$
L^{'}(z)\cdot (L(z)-a)^{N}=e^{p(z)}\cdot e^{-G(L(z))},
$$
$$
\log L^{'}(z)+N\log(L(z)-a)=p(z)-G(L(z)).
$$
As in the proof of Theorem \ref{thm6}, this can not have a solution in which $G(z)\not\equiv {\rm Const.}$.
So $G(z)\equiv {\rm Const.}$ (it can even be equal to $0$) and hence necessarily $N\in\mathbb{Z}^{+}$. Our equation
becomes:
$$
L^{'}(z)\cdot (L(z)-a)^{N}=\omega_{0}\cdot e^{p(z)},\,\,\,\,\,\omega_{0}\in\mathbb{C}^{\times}.
$$
Hence by integrating we get:
\begin{equation}
\label{eq13}
\left(\frac{1}{N+1}\right)\cdot (L(z)-a)^{N+1}=f(z).
\end{equation}
This has a solution which is an entire function if and only if $f(z)\ne 0$ $\forall\,z\in\mathbb{C}$, i.e. $f(z)$
has a Picard exceptional value (by our normalization it is $a=0$). In this case:
\begin{equation}
\label{eq14}
L(z)=(N+1)^{1/(N+1)}\cdot f(z)^{1/(N+1)}+a,
\end{equation}
is indeed in ${\rm elh}(\mathbb{C})-{\rm Aut}(\mathbb{C})$ (by differentiating equation (\ref{eq14}) 
and $L(\mathbb{C})=\mathbb{C}-\{ a\}$. Equation (\ref{eq13}) gives a different factorization (of $f(z)$), then the
one given in Theorem \ref{thm4}. Namely:
\begin{equation}
\label{eq15}
f(z)=\left(\frac{z}{N+1}\right)\circ z^{N+1}\circ (z-a)\circ L(z).
\end{equation}
Theorem \ref{thm7} is now proved. $\qed $

\begin{corollary}\label{cor1}
The functions $f_{2k}(z)=\int_{0}^{z}e^{-t^{2k}}dt$, $k\in\mathbb{Z}^{+}$, belong to ${\rm elh}(\mathbb{C})-{\rm Aut}(\mathbb{C})$,
satisfy $f_{2k}(\mathbb{C})=\mathbb{C}$, and are prime entire mappings.
\end{corollary}
\noindent
We recall that after Definition \ref{def4} and within Remark \ref{rem1} we denoted by $\mathcal{M}$ the set of all the 
(equivalence classes of the) primes in ${\rm elh}(\mathbb{C})$. In Theorem \ref{thm1} we gave the two (right and left)
fractal representations of ${\rm elh}(\mathbb{C})$. Namely:
$$
{\rm elh}(\mathbb{C})=W_{R}\cup\bigcup_{p\in\mathcal{M}}R_{p}({\rm elh}(\mathbb{C}))=
W_{L}\cup\bigcup_{p\in\mathcal{M}}L_{p}({\rm elh}(\mathbb{C})).
$$
We have seen (in Proposition \ref{prop6}) that $e^{z}\in W_{R}-({\rm Aut}(\mathbb{C})\cup\mathcal{M})$, and also 
$e^{z}\in W_{L}-({\rm Aut}(\mathbb{C})\cup\mathcal{M})$. Remark \ref{rem4} suggested the unpleasant possibility
that the fractal representations of Theorem \ref{thm1} might degenerate to:
$$
{\rm elh}(\mathbb{C})=W_{R}\,\,\,\,\,{\rm or}\,\,\,\,\,{\rm elh}(\mathbb{C})=W_{L}.
$$
This happens exactly when ${\rm elh}(\mathbb{C})$ contains no primes, i.e. when $\mathcal{M}=\emptyset$. In fact
the investigation of this possibility was the author's main motivation to conduct the research in this manuscript.
We can now finally deduce that these two non desirable situations can not occur:

\begin{theorem}\label{thm8}
\begin{equation}
\label{eq16}
\mathcal{M}\ne\emptyset.
\end{equation}
\end{theorem}
\noindent
{\bf Proof.} \\
This follows by Corollary \ref{cor1}. $\qed $ \\

\noindent
{\it Ronen Peretz \\
Department of Mathematics \\ Ben Gurion University of the Negev \\
Beer-Sheva , 84105 \\ Israel \\ E-mail: ronenp@math.bgu.ac.il} \\ 
 
\end{document}